\pdfoutput=1
\RequirePackage{ifpdf}
\ifpdf
\documentclass[pdftex]{sigma}
\else
\documentclass{sigma}
\fi

\usepackage{aligned-overset}

\numberwithin{equation}{section}

\newtheorem{Theorem}{Theorem}[section]
\newtheorem*{Theorem*}{Theorem}
\newtheorem{Corollary}[Theorem]{Corollary}
\newtheorem{Lemma}[Theorem]{Lemma}
\newtheorem{Proposition}[Theorem]{Proposition}

\theoremstyle{definition}
\newtheorem{Definition}[Theorem]{Definition}
\newtheorem*{Notation}{Notation}
\newtheorem{Example}[Theorem]{Example}
\newtheorem{Remark}[Theorem]{Remark}

\def\id{\operatorname{id}}

\def\Aut{\operatorname{Aut}}
\def\Inn{\operatorname{Inn}}
\def\Der{\operatorname{Der}}
\def\exp{\operatorname{exp}}
\def\d{\mathrm{d}}

\def\Hom{\operatorname{Hom}}

\def\End{\operatorname{End}}

\newcommand{\nc}{\newcommand}
\nc{\tred}[1]{\textcolor{red}{#1}}
\nc{\tblue}[1]{\textcolor{blue}{#1}}
\nc{\tgreen}[1]{\textcolor{green}{#1}}
\nc{\tpurple}[1]{\textcolor{purple}{#1}}
\nc{\btred}[1]{\textcolor{red}{\bf #1}}
\nc{\btblue}[1]{\textcolor{blue}{\bf #1}}
\nc{\btgreen}[1]{\textcolor{green}{\bf #1}}
\nc{\btpurple}[1]{\textcolor{purple}{\bf #1}}
\nc{\NN}{{\mathbb N}}
\nc{\ncsha}{{\mbox{\cyr X}^{\mathrm NC}}} \nc{\ncshao}{{\mbox{\cyr
X}^{\mathrm NC}_0}}

\renewcommand{\frak}{\mathfrak}
\newcommand{\delete}[1]{}

\nc{\mlabel}[1]{\label{#1}}
\nc{\mcite}[1]{\cite{#1}}
\nc{\mref}[1]{\ref{#1}}
\nc{\meqref}[1]{\eqref{#1}}
\nc{\mbibitem}[1]{\bibitem{#1}}
\delete{% Use the next lines to show names
\nc{\mlabel}[1]{\label{#1}{\hfill \hspace{1cm}{\bf{{\ }\hfill(#1)}}}}
\nc{\mcite}[1]{\cite{#1}{{\bf{{\ }(#1)}}}}
\nc{\mref}[1]{\ref{#1}{{\bf{{\ }(#1)}}}}
\nc{\meqref}[1]{\eqref{#1}{{\bf{{\ }(#1)}}}}
\nc{\mbibitem}[1]{\bibitem[\bf #1]{#1}}
}
\nc{\Shu}{\mathrm{Sh}}
\nc{\tot}{\big|}
\nc{\AAO}{{A,\Omega}}
\nc{\msha}{\sha_\Omega}
\nc{\mrsha}{\sha_{\Omega,F}}

\nc{\opa}{\ast} \nc{\opb}{\odot} \nc{\op}{\bullet} \nc{\pa}{\frakL}
\nc{\arr}{\rightarrow} \nc{\lu}[1]{(#1)} \nc{\mult}{\mrm{mult}}
\nc{\diff}{\mathfrak{Diff}}
\nc{\opc}{\sharp}\nc{\opd}{\natural}
\nc{\ope}{\circ}
\nc{\dpt}{\mathrm{d}}
\nc{\hck}{H_{RT}}
\nc{\vdf}{\calf}
\nc{\ldf}{\calf_\ell}
\nc{\hlf}{H_\ell}
\nc{\onek}{\mathbf{1}_\bfk}
\nc{\tforall}{\, \text{ for }\, }
\nc{\qforall}{\quad \text{for all }}

\nc{\mrba}{MRBA\xspace}
\nc{\Mrba}{MRBA\xspace}
\nc{\mrbas}{MRBAs\xspace}
\nc{\Mrbas}{MRBAs\xspace}
\nc{\match}{matching\xspace}
\nc{\Match}{Matching\xspace}
\nc{\Mza}{Matching Zinbiel algebra\xspace}
\nc{\Mzas}{Matching Zinbiel algebras\xspace}
\nc{\mza}{matching Zinbiel algebra\xspace}
\nc{\mzas}{matching Zinbiel algebras\xspace}

\nc{\za}{Zinbiel algebra\xspace}
\nc{\paybe}{polarized associative Yang--Baxter equation\xspace}
\nc{\Paybe}{Polarized associative Yang--Baxter equation\xspace}
\nc{\cpaybe}{PAYBE}

\nc{\diam}{alternating\xspace}
\nc{\Diam}{Alternating\xspace}
\nc{\cdiam}{canonical alternating\xspace}
\nc{\Cdiam}{Canonical alternating\xspace}
\nc{\AW}{\mathcal{A}}
\nc{\rba}{Rota--Baxter algebra\xspace}
\nc{\rbas}{Rota--Baxter algebras\xspace}

\nc{\ari}{\mathrm{ar}}
\nc{\lef}{\mathrm{lef}}
\nc{\Sh}{\mathrm{ST}}
\nc{\Cr}{\mathrm{Cr}}
\nc{\st}{{Schr\"oder tree}\xspace}
\nc{\sts}{{Schr\"oder trees}\xspace}
\nc{\vertset}{\Omega} % set of vertex decorations
\nc{\assop}{\quad \begin{picture}(5,5)(0,0)
\line(-1,1){10}
\put(-2.2,-2.2){$\bullet$}
\line(0,-1){10}\line(1,1){10}
\end{picture} \quad \smallskip}

\nc{\operator}{\begin{picture}(5,5)(0,0)
\line(0,-1){6}
\put(-2.6,-1.8){$\bullet$}
\line(0,1){9}
\end{picture}}

\nc{\idx}{\begin{picture}(6,6)(-3,-3)
\put(0,0){\line(0,1){6}}
\put(0,0){\line(0,-1){6}}
\end{picture}}

\nc{\pb}{{\mathrm{pb}}}
\nc{\Lf}{{\mathrm{Lf}}}

\nc{\lft}{{left tree}\xspace}
\nc{\lfts}{{left trees}\xspace}

\nc{\fat}{{fundamental averaging tree}\xspace}

\nc{\fats}{{fundamental averaging trees}\xspace}
\nc{\avt}{\mathrm{Avt}}

\nc{\rass}{{\mathit{RAss}}}

\nc{\aass}{{\mathit{AAss}}}

\nc{\twovec}[2]{{\textstyle \left[#1\atop #2\right]}}

\nc{\vin}{{\mathrm Vin}} %decoration set of indices
\nc{\lin}{{\mathrm Lin}} %decoration set of leaves
\nc{\inv}{\mathrm{I}n}
\nc{\gensp}{V} % space of generators
\nc{\genbas}{\mathcal{V}} % basis of the space of generators
\nc{\bvp}{V_P} % Rota--Baxter generating space
\nc{\gop}{{\,\omega\,}} % generic binary operation

\nc{\bin}[2]{ (_{\stackrel{\scs{#1}}{\scs{#2}}})} %binomial coeff
\nc{\binc}[2]{ \left (\!\! \begin{array}{c} \scs{#1}\\
 \scs{#2} \end{array}\!\! \right )} %binomial coeff
\nc{\bincc}[2]{ \left ( {\scs{#1} \atop
 \vspace{-1cm}\scs{#2}} \right )} %binomial coeff
\nc{\bs}{\bar{S}} \nc{\cosum}{\sqsubset} \nc{\la}{\longrightarrow}
\nc{\rar}{\rightarrow} \nc{\dar}{\downarrow} \nc{\dprod}{**}
\nc{\dap}[1]{\downarrow \rlap{$\scriptstyle{#1}$}}
\nc{\md}{\mathrm{dth}} \nc{\uap}[1]{\uparrow
\rlap{$\scriptstyle{#1}$}} \nc{\defeq}{\stackrel{\rm def}{=}}
\nc{\disp}[1]{\displaystyle{#1}} \nc{\dotcup}{\
\displaystyle{\bigcup^\bullet}\ } \nc{\gzeta}{\bar{\zeta}}
\nc{\hcm}{\ \hat{,}\ } \nc{\hts}{\hat{\otimes}}
\nc{\barot}{{\otimes}} \nc{\free}[1]{\bar{#1}}
\nc{\uni}[1]{\tilde{#1}} \nc{\hcirc}{\hat{\circ}} \nc{\lleft}{[}
\nc{\lright}{]} \nc{\lc}{\lfloor} \nc{\rc}{\rfloor}
\nc{\curlyl}{\left \{ \begin{array}{c} {} \\ {} \end{array}
 \right . \!\!\!\!\!\!\!}
\nc{\curlyr}{ \!\!\!\!\!\!\!
 \left . \begin{array}{c} {} \\ {} \end{array}
 \right \} }
\nc{\longmid}{\left | \begin{array}{c} {} \\ {} \end{array}
 \right . \!\!\!\!\!\!\!}
\nc{\onetree}{\bullet} \nc{\ora}[1]{\stackrel{#1}{\rar}}
\nc{\ola}[1]{\stackrel{#1}{\la}}%${\Bbb Z}$
\nc{\ot}{\otimes} \nc{\mot}{{{\boxtimes\,}}}
\nc{\otm}{\overline{\boxtimes}} \nc{\sprod}{\bullet}
\nc{\scs}[1]{\scriptstyle{#1}} \nc{\mrm}[1]{{\rm #1}}
\nc{\margin}[1]{\marginpar{\rm #1}} %{\rm #1}}
\nc{\dirlim}{\displaystyle{\lim_{\longrightarrow}}\,}
\nc{\invlim}{\displaystyle{\lim_{\longleftarrow}}\,}
\nc{\mvp}{\vspace{0.3cm}} \nc{\tk}{^{(k)}} \nc{\tp}{^\prime}
\nc{\ttp}{^{\prime\prime}} \nc{\svp}{\vspace{2cm}}
\nc{\vp}{\vspace{8cm}} \nc{\proofbegin}{\noindent{\bf Proof: }}
\nc{\proofend}{$\blacksquare$ \vspace{0.3cm}}
\nc{\modg}[1]{\!<\!\!{#1}\!\!>}
\nc{\intg}[1]{F_C(#1)} \nc{\lmodg}{\!
<\!\!} \nc{\rmodg}{\!\!>\!}
\nc{\cpi}{\widehat{\Pi}}
\nc{\sha}{{\mbox{\cyr X}}} %used to be \cyr

\newfont{\scyr}{wncyr10 scaled 550}
\nc{\ssha}{\mbox{\bf \scyr X}}
\nc{\shap}{\,{\mbox{\cyrs X}}\,} %sha as product
\nc{\shpr}{\diamond} %Shuffle product
\nc{\shp}{\ast} \nc{\shplus}{\shpr^+}
\nc{\shprc}{\shpr_c} %Cartier's product
\nc{\msh}{\ast} \nc{\zprod}{m_0} \nc{\oprod}{m_1}
\nc{\vep}{\epsilon} \nc{\labs}{\mid\!} \nc{\rabs}{\!\mid}
\nc{\sqmon}[1]{\langle #1\rangle}
\nc{\mmbox}[1]{\mbox{\ #1\ }} \nc{\dep}{\mrm{dep}} \nc{\fp}{\mrm{FP}}
\nc{\Mor}{Mor\xspace} \nc{\gmzvs}{gMZV\xspace}
\nc{\gmzv}{gMZV\xspace} \nc{\mzv}{MZV\xspace}
\nc{\im}{\mrm{im}} \nc{\incl}{\mrm{incl}} \nc{\map}{\mrm{Map}}
\nc{\mchar}{\rm char} \nc{\nz}{\rm NZ} \nc{\supp}{\mathrm Supp}

\nc{\Alg}{\mathbf{Alg}} \nc{\Bax}{\mathbf{Bax}} \nc{\bff}{\mathbf f}
\nc{\bfk}{{\bf k}} \nc{\bfone}{{\bf 1}} \nc{\bfx}{\mathbf x}
\nc{\bfy}{\mathbf y}
\nc{\base}[1]{\bfone^{\otimes ({#1}+1)}} %{{a_{#1}}}
\nc{\Cat}{\mathbf{Cat}}

\nc{\detail}{\marginpar{\bf More detail}
 \noindent{\bf Need more detail!}
 \svp}
\nc{\Int}{\mathbf{Int}} \nc{\Mon}{\mathbf{Mon}}
\nc{\rbtm}{{shuffle }} \nc{\rbto}{{Rota--Baxter }}
\nc{\remarks}{\noindent{\bf Remarks: }} \nc{\Rings}{\mathbf{Rings}}
\nc{\Sets}{\mathbf{Sets}} \nc{\wtot}{\widetilde{\odot}}
\nc{\wast}{\widetilde{\ast}} \nc{\bodot}{\bar{\odot}}
\nc{\bast}{\bar{\ast}} \nc{\hodot}[1]{\odot^{#1}}
\nc{\hast}[1]{\ast^{#1}} \nc{\mal}{\mathcal{O}}
\nc{\tet}{\tilde{\ast}} \nc{\teot}{\tilde{\odot}}
\nc{\oex}{\overline{x}} \nc{\oey}{\overline{y}}
\nc{\oez}{\overline{z}} \nc{\oef}{\overline{f}}
\nc{\oea}{\overline{a}} \nc{\oeb}{\overline{b}}
\nc{\weast}[1]{\widetilde{\ast}^{#1}}
\nc{\weodot}[1]{\widetilde{\odot}^{#1}} \nc{\hstar}[1]{\star^{#1}}
\nc{\lae}{\langle} \nc{\rae}{\rangle}
\nc{\lf}{\lfloor}
\nc{\rf}{\rfloor}
\nc{\QQ}{{\mathbb Q}}
\nc{\RR}{{\mathbb R}} \nc{\ZZ}{{\mathbb Z}}

\nc{\cala}{{\mathcal A}} \nc{\calb}{{\mathcal B}}
\nc{\calc}{{\mathcal C}}
\nc{\cald}{{\mathcal D}} \nc{\cale}{{\mathcal E}}
\nc{\calf}{{\mathcal F}} \nc{\calg}{{\mathcal G}}
\nc{\calh}{{\mathcal H}} \nc{\cali}{{\mathcal I}}
\nc{\call}{{\mathcal L}} \nc{\calm}{{\mathcal M}}
\nc{\caln}{{\mathcal N}} \nc{\calo}{{\mathcal O}}
\nc{\calp}{{\mathcal P}} \nc{\calr}{{\mathcal R}}
\nc{\cals}{{\mathcal S}} \nc{\calt}{{\mathcal T}}
\nc{\calu}{{\mathcal U}} \nc{\calw}{{\mathcal W}} \nc{\calk}{{\mathcal K}}
\nc{\calx}{{\mathcal X}} \nc{\CA}{\mathcal{A}}

\nc{\fraka}{{\mathfrak a}} \nc{\frakA}{{\mathfrak A}}
\nc{\frakb}{{\mathfrak b}} \nc{\frakB}{{\mathfrak B}}
\nc{\frakc}{{\mathfrak c}}
\nc{\frakD}{{\mathfrak D}} \nc{\frakF}{\mathfrak{F}}
\nc{\frakf}{{\mathfrak f}} \nc{\frakg}{{\mathfrak g}}
\nc{\frakH}{{\mathfrak H}} \nc{\frakL}{{\mathfrak L}}
\nc{\frakM}{{\mathfrak M}} \nc{\bfrakM}{\overline{\frakM}}
\nc{\frakm}{{\mathfrak m}} \nc{\frakP}{{\mathfrak P}}
\nc{\frakN}{{\mathfrak N}} \nc{\frakp}{{\mathfrak p}}
\nc{\frakS}{{\mathfrak S}} \nc{\frakT}{\mathfrak{T}}
\nc{\frakX}{{\mathfrak X}}
\nc{\BS}{\mathbb{S
}}

\font\cyr=wncyr10 \font\cyrs=wncyr7
\nc{\li}[1]{\textcolor{red}{#1}}
\nc{\lir}[1]{\textcolor{red}{Li:#1}}
\nc{\yi}[1]{\textcolor{blue}{Yi: #1}}
\nc{\xing}[1]{\textcolor{purple}{Xing:#1}}
\nc{\revise}[1]{\textcolor{red}{#1}}

\nc{\ID}{{\rm I}}\nc{\lbar}[1]{\overline{#1}}\nc{\bre}{{\rm bre}}
\nc{\sd}{\cals}\nc{\rb}{\rm RB}\nc{\A}{\rm A}\nc{\LL}{\rm L}\nc{\tx}{\tilde{X}}
\nc{\col}{\Delta_{RT}}\nc{\mul}{m_{RT}}\nc{\ul}{u_{RT}}\nc{\epl}{\varepsilon_{RT}}
\nc{\hl}{H_{RT}}\nc{\arro}[1]{#1}\nc{\px}{P_{\tx}}\nc{\pw}{P_{\mathfrak{w}}}\nc{\pl}{B_\omega^+}
\nc{\pp}{\pl}\nc{\ppp}[1]{B^+(#1)}\nc{\dw}{\diamond_{\mathfrak{w}}}\nc{\dl}{\diamond_{\rm \ell}}
\nc{\ncshaw}{\sha^{{\rm NC}}_{\Omega}}\nc{\ncshal}{\sha^{{\rm NC}}_{{\rm RT}}}
\nc{\ver}{\rm V}\nc{\ld}{l}\nc{\del}{\Delta_{{\rm \ell}}}\nc{\epsl}{\epsilon_{{\rm \ell}}}
\nc{\uul}{u_{{\rm \ell}}}\nc{\oneh}{\mathbf{1}}\nc{\onew}{\mathbf{1}}
\nc{\etree}{1} \nc{\conc}{m_{RT}}
\nc{\hrtb}{H_{RT}(X\sqcup\Omega)} \nc{\hrts}{H_{RT}(X, \Omega)}\nc{\rts}{\mathcal{T}(X, \Omega)}\nc{\rfs}{\mathcal{F}(X, \Omega)} \nc{\ncshall}{\sha^{{\rm NC}}_{{\rm RT}}} \nc{\ldl}{\leq_{\mathrm{dl}}} \nc{\pla}{B_{\alpha}^{+}} \nc{\plb}{B_{\beta}^{+}}

\nc{\bim}[1]{#1} \nc{\shaop}{\sha_{\Omega}^{+}} \nc{\shao}{\sha_{\Omega}}
\nc{\bbim}[2]{#1 #2} \nc{\bbbim}[2]{#1,\, #2} \nc{\RBF}{{\rm MRB}}
\nc{\frbf}{F_{\RBF}} \nc{\shaf}{\ssha_{\tiny{\Omega}}} \nc{\sham}{\diamond}
\nc{\dnx}{\Delta_n A} \nc{\dx}{\Delta A} \nc{\dgp}{{\rm deg_{P}}}
\nc{\dgt}{{\rm deg_{T}}} \nc{\dg}{{\rm deg}} \nc{\ida}{ID($A$)} \nc{\tu}{\tilde{u}} \nc{\tv}{\tilde{v}}
 \nc{\fbase}{\calb} \nc{\LF}{\mathrm{RF}} \nc{\FFA}{\mathrm{LF}} \nc{\irr}{\mathrm{Irr}}
 \nc{\result}{\bfk\mathrm{Irr}(S_n)} \nc{\I}{I_{\mathrm{ID},n}^0}
 \nc{\nrs}{\calr_n^\star} \nc{\ii}{\mathrm{I}} \nc{\iii}{\mathrm{II}}
\nc{\intl}{{\rm int}}\nc{\ws}[1]{{#1}}\nc{\deleted}[1]{\delete{#1}}\nc{\plas}{placements\xspace}

\nc{\Id}{\mathrm{Id}} \nc{\Irr}{\mathrm{Irr}}

\nc{\tos}{totally ordered set} \nc{\nes}{nonempty set}
\nc{\rsha}{\sha^{\rm rel}}
\nc{\lra}{\longrightarrow}

\nc{\rel}{\mathrm{\text{rel}}}

\nc{\rsham}{\sham}
\nc{\bo}{\tau}

\nc{\sshao}{~\ssha_{\Omega}~}
\nc{\resham}{\diamond_{\Omega}^{\mathrm{rel}}}

\nc{\fraku}{\frak{u}}
\nc{\frakv}{\frak{v}}
\nc{\frabu}{\bar{\mathfrak{u}}}
\nc{\frabv}{\bar{\mathfrak{v}}}
\nc{\frabw}{\bar{\mathfrak{w}}}

\nc{\udl}[1]{\underline{#1}}

\nc{\Po}{(P_\omega)_{\omega\in \Omega}}
\nc{\basf}{F} \nc{\frakhat}[1]{\frak{#1'}}
\nc{\cten}[2]{\left[\begin{array}{c}#1 \\#2 \\ \end{array}\right]}

\renewcommand{\d}{\operatorname{d}}

\begin{document}

\allowdisplaybreaks

\newcommand{\arXivNumber}{2307.10535}

\renewcommand{\PaperNumber}{099}

\FirstPageHeading

\ShortArticleName{Twisted Post-Groups and Skew Trusses}

\ArticleName{Twisted Post-Groups and Skew Trusses}

\Author{Shukun WANG~$^{\rm ab}$}

\AuthorNameForHeading{S.~Wang}

\Address{$^{\rm a)}$~School of Mathematics and Big Data, Anhui University of Science and Technology,\\
\hphantom{$^{\rm a)}$}~Huainan 232001, P.R.~China}
\EmailD{\mail{2024093@aust.edu.cn}}
\Address{$^{\rm b)}$~Anhui Province Engineering Laboratory for Big Data Analysis and Early Warning Technology\\
\hphantom{$^{\rm b)}$}~of Coal Mine Safety, Huainan 232001, P.R.~China}

\ArticleDates{Received July 13, 2026, in final form October 09, 2026; Published online October 19, 2026}

\Abstract{To understand the origin of post-groups introduced by Bai, Guo, Sheng, and Tang from the perspective of rings, we introduce the notion of (weak) twisted post-groups. First, we show that every element of a twisted post-group belongs to a unique group, and that a twisted post-group can be viewed as the disjoint union of such groups. Next, we prove that the category of weak twisted post-groups is isomorphic to the category of skew trusses, and that every two-sided twisted post-group admits the structure of a two-sided skew brace. It~follows that every abelian two-sided twisted post-group gives rise to a radical ring. We~then introduce twisted post-Lie algebras and investigate their algebraic properties, showing that differentiating a twisted post-Lie group yields a twisted post-Lie algebra. Finally, we~consider the linearization of (weak) twisted post-groups and propose the notion of (weak) twisted post-Hopf algebras. We~show that every twisted post-Hopf algebra gives rise to another Hopf algebra, called the sub-adjacent Hopf algebra.}

\Keywords{(weak) twisted post-group; skew truss; Rota--Baxter system of groups; twisted post-Lie algebra}

\Classification{22E60; 17B40; 16T05}

\vspace{-2mm}

\section{Introduction}\label{Sec:Intro}
\subsection{Braces and trusses}
In modern theoretical and mathematical physics, the Yang--Baxter equation plays a fundamental role. The Yang--Baxter equation has its origins in the study of exactly solvable models in statistical mechanics \cite{Baxter1} and in Yang’s investigation of interacting particle systems~\cite{Yang}. The investigation into the set-theoretical solutions of the Yang--Baxter equation was pioneered by Drinfeld~\cite{Drinfeld} and pursued by several authors \cite{Etingof,TG,LJ}.

In \cite{Rump1}, Rump introduced braces as a generalisation of radical rings, which yield non-degen\-er\-ate involutive set-theoretic solutions of the Yang--Baxter equation. Subsequently, in~\cite{Cedo}, Ced\'o, Jespers and Okni\'nski defined left braces. For further developments on braces and their connection to set-theoretic solutions of the Yang--Baxter equation, we refer to~\cite{Cedo1,Rump1,SM}. More recently, Guarnieri and Vendramin introduced the notion of skew braces, which generalise braces to the non-abelian setting. Later, Bachiller proved that every non-degenerate set-theoretic solution of the Yang--Baxter equation can be constructed from skew braces~\cite{BA}.

To understand how two group operations interact in a skew brace, Brzezi\'nski~\cite{TB2} proposed a~new algebraic system, called a skew truss. The notion of a Rota--Baxter system of groups was introduced as the skew truss analogue of Rota--Baxter groups~\cite{ZH}.

More recently, the notion of a skew truss was also introduced in a broader
categorical framework in~\cite{RG}, where Hopf trusses were
introduced in an arbitrary braided monoidal category~$\mathcal{C}$.
When $\mathcal{C}$ is the category of sets, Hopf trusses are just skew trusses~\cite{RG}.

\subsection{Post-Lie algebras and post-groups}

Post-type algebraic structures appear in various areas of mathematics and physics.
In~\cite{VB}, Vallette introduced the notion of post-Lie algebras in the study of operads. In \cite{BC}, Bai, Guo and Ni showed that post-Lie algebras are closely related to the classical Yang--Baxter equation. Then in \cite{KA}, Munthe-Kaas and Lundervold investigated the connections of post-Lie algebras with Lie-Butcher series and flows on manifolds. For more works on post-Lie algebras, we refer to \cite{YB,DB,DV}.

Recently, in \cite{YNL}, Li, Sheng and Tang introduced the Hopf algebra analogues of post-Lie algebras, called post-Hopf algebras. It was shown that the category of post-Hopf algebras and the category of Hopf braces introduced in \cite{IA} are isomorphic under cocommutativity conditions. More recently in~\cite{BM}, Bai, Guo, Sheng and Tang introduced the notion of post-groups, which is the group analogue of post-Lie algebras. A post-group is a triple $(G,\cdot,\rhd)$, where $(G,\cdot)$ is a~group and $\rhd\colon G\times G\to G$ is an operation such that
\begin{enumerate}\itemsep=0pt
	\item For each $a\in G$, the left multiplication $L^{\rhd}_{a}\colon G\to G$ given by $L^{\rhd}_a(b)=a\rhd b$, $\forall b\in G,$ is an automorphism;
	\item For any $a,b,c\in G$, $(a\cdot(a\rhd b))\rhd c=a\rhd (b\rhd c)$.
\end{enumerate}
It was shown that Butcher groups in numerical integration, $\mathcal{P}$-groups of operads, and braided groups all have the structure of post-groups. Moreover, it was proven that the category of post-groups is isomorphic to the category of skew braces.

\subsection{The skew truss analogue of post-groups}

In view of the equivalence between the category of post-groups and the category of skew braces, and the fact that every skew brace admits a skew truss structure, we introduce in this paper the notion of (weak) twisted post-groups, defined as the skew truss analogue of post-groups.

First, we give some examples of (weak) twisted post-groups and show that in any twisted post-group $(G,\cdot,\rhd,\Phi)$, there exists an operation, called the sub-adjacent operation, under which~$G$ forms a semigroup.

Then we prove that for any $a$ in a twisted post-group, there is a~unique group~$G_a$ such that~$G_a$ is a group with respect to the sub-adjacent operation. Furthermore, we give a decomposition theorem for twisted post-groups, which generalises the decomposition theorem of Rota--Baxter systems of groups.

Next, we investigate the relationship between weak twisted post-groups and skew trusses. We prove that the category of weak twisted post-groups and the category of skew trusses are isomorphic. Furthermore, we show that a skew truss carries the structure of a twisted post-group if and only if it is right divisible with respect to its sub-adjacent operation. Then we investigate the connections between twisted post-groups and Rota--Baxter systems of groups. We also study the relationship between (weak) twisted post-groups and rings. In particular, every abelian two-sided twisted post-group gives rise to a radical ring.

As the Lie algebra analogue of twisted post-groups, we introduce the notion of twisted post-Lie algebras. Twisted post-Lie algebras can be regarded as the skew truss version of post-Lie algebras. We prove that differentiating a twisted post-Lie group yields a twisted post-Lie algebra. Finally, as the Hopf algebra analogue of twisted post-groups, we introduce (weak) twisted post-Hopf algebras. We show that the category of weak twisted post-Hopf algebras is isomorphic to the category of Hopf trusses, and that every twisted post-Hopf algebra gives rise to another Hopf algebra, called its sub-adjacent Hopf algebra.

\subsection{Outline of this paper}
This paper is organised as follows. In Section~\ref{2}, we introduce (weak) twisted post-groups and provide several examples. Section~\ref{3} is devoted to the study of their internal structure. In Section~\ref{4}, we investigate the connections between skew trusses and (weak) twisted post-groups. The relationship between post-groups and rings is also discussed. In Section~\ref{5}, we define twisted post-Lie algebras and consider their connections with the differentiations of twisted post-Lie groups. In Section~\ref{6}, we introduce (weak) twisted post-Hopf algebras as the Hopf algebra analogues of (weak) twisted post-Lie algebras. Finally, we show that the category of weak twisted post-Hopf algebras is isomorphic to the category of Hopf trusses.

\begin{Notation} Throughout this paper, for a group $(G,\cdot)$, denote its identity by $1$. Given a unital associative algebra $A$ over a field $\mathbb{F}$, we also denote its identity by $1$. And the unit map $\mathbb{F}\to A$ given by $a\to a1$ is also denoted by $1$. For any coalgebra $(C,\Delta,\epsilon)$, we will use the Sweedler
notation $\Delta(a)=a_{1}\otimes a_{2}$ without summation sign. Thanks to the coassociativity, we can extend the notation to write
\[a_1\otimes a_2\otimes a_3=a_{1,1}\otimes a_{1,2}\otimes a_2=a_1\otimes a_{2,1}\otimes a_{2,2}.\] Similarly for any $n
>3$, we can write
\[a_1\otimes a_2\otimes\cdots\otimes a_n=\bigl(\Delta\otimes \id^{n-1}\bigr)\cdots(\Delta\otimes \id)\Delta(a)=\bigl( \id^{n-1}\otimes \Delta\bigr)\cdots(\id\otimes\Delta)\Delta(a).\]
In this notation, the counitality can be written as
\[\epsilon(a_1)a_2=a=a_1\epsilon(a_2),\qquad\forall a\in C.\]
Unless otherwise cited, a Hopf algebra $(H,\cdot,\Delta,1,\epsilon,S)$ is simply denoted by $H$.
\end{Notation}

\section{(Weak) Twisted post-groups}\label{2}
In this section, we introduce the notion of (weak) twisted post-groups as the skew truss analogue of post-groups. Several examples are provided, and some fundamental algebraic properties of twisted post-groups are investigated.

Generalising the notion of post-groups introduced in~\cite{BM}, we propose the notion of (weak) twisted post-groups.
\begin{Definition}\label{TPG}
	Let $(G,\cdot)$ be a group with operators $\rhd\colon G\times G\to G$ and $\Phi\colon G\to G$. Define $\circ\colon G\times G\to G$ by \[a\circ b=\Phi(a)\cdot(a\rhd b),\qquad\forall a,b\in G.\] $(G,\cdot,\rhd,\Phi)$ is called a \textit{left twisted post-group} if for any $a,b,c\in G$,
	\begin{enumerate}\itemsep=0pt
		\item the left multiplication $L_{a}^{\rhd}$ defined by
		\begin{equation}\label{TPG1}
			L_{a}^{\rhd}(b)=a\rhd b,
		\end{equation}
		is a group automorphism of $(G,\cdot)$;
		\item the following ``twisted weighted associativity'' holds:
		\begin{equation}\label{TPG2}
			(a\circ b)\rhd c=a\rhd(b\rhd c);
		\end{equation}
		\item the ``left compatible condition for $\Phi$'' holds:
		\begin{equation}\label{TPG3}
			\Phi(a\circ b)=a\circ \Phi(b).
		\end{equation}
	\end{enumerate}
	The operator $\Phi$ is called the \textit{cocycle} and $\circ$ is called the \textit{sub-adjacent operation} of $(G,\cdot,\rhd,\Phi)$. If condition (a) is weakened by requiring $L_a^\rhd$ to be a group endomorphism, not necessarily an automorphism, then $(G,\cdot,\rhd,\Phi)$ is called a \textit{left weak twisted post-group}. Note that if $\Phi$ is the identity map of $G$, then the left twisted post-group is a post-group.
\end{Definition}

In a symmetric way, let $(G,\cdot)$ be a group with operators $\lhd\colon G\times G\to G$ and $\Phi\colon G\to G$. Define $\circ\colon G\times G\to G$ by
\[a\circ b=(a\lhd b)\cdot\Phi(b),\qquad \forall a,b\in G.\]
The quadruple $(G,\cdot,\lhd,\Phi)$ is called a right twisted post-group if
\begin{enumerate}\itemsep=0pt
	\item the right multiplication $R_{b}^{\lhd}$ defined by $R_{b}^{\lhd}(a)=a\lhd b$,
	is a group automorphism of $(G,\cdot)$;
	\item the following ``twisted weighted associativity'' holds:
	$a\lhd(b\circ c)=(a\lhd b)\lhd c$;
	\item the ``right compatible condition for $\Phi$'' holds:
$\Phi(a\circ b)=\Phi(a)\circ b$,
\end{enumerate}
for any $a,b,c\in G$. $(G,\cdot,\lhd,\Phi)$ is called a \textit{right weak twisted post-group} if $R^{\lhd}_{b}$ is an endomorphism of $(G,\cdot)$.

Without further mention, we will simply refer to left (weak) twisted post-groups as (weak) twisted post-groups in this note.

In the next example, we show that every Rota--Baxter system of groups introduced in \cite{ZH} has the structure of a twisted post-group.

\begin{Example}\label{RG}
	Let $(G,\cdot)$ be a group and $(G,B_1,B_2)$ be a Rota--Baxter system of groups, that is, $B_1\colon G\to G$ and $B_2\colon G\to G$ satisfy
	\[B_1(a)B_1(b)=B_1(B_1(a)bB_2(a)),\qquad
 B_2(b)B_2(a)=B_2(B_1(a)bB_2(a)),\] for any $a,b\in G$. Define $\rhd\colon G\times G\to G$ by \[a\rhd b=B_2(a)^{-1}bB_2(a),\qquad \forall a,b\in G\] and $\Phi\colon G\to G$ by
	\[\Phi(a)=B_1(a)B_2(a),\qquad\forall a\in G.\] Then $(G,\cdot,\rhd,\Phi)$ is a twisted post-group. 	
	Similarly, one can define
	$\lhd\colon G\times G\to G$ by \[a\lhd b=B_1(b)aB_1(b)^{-1},\qquad \forall a,b\in G,\] and $\Phi'\colon G\to G$ by
	\[\Phi'(a)=B_1(a)B_2(a),\qquad\forall a\in G,\] which yields a right twisted post-group $(G,\cdot,\lhd,\Phi')$. 	
\end{Example}

In~\cite{BM}, it was shown that every operad is associated with a post-group. The next example shows that every operad also gives rise to a weak twisted post-group with trivial cocycle.
\begin{Example}
	An operad $\{ \mathcal{P}(n)\}_{n\geq 1}$ is a sequence of $\mathbb{S}_n$-modules with a family of linear operators
	\[\gamma\colon \ \mathcal{P}(n)\otimes\mathcal{P}(k_1)\otimes\cdots\otimes\mathcal{P}(k_n)\to \mathcal{P}(k_1+k_2+\cdots+k_n),\qquad\forall n,k_1,\dots,k_n\ge 1,\]
	such that they satisfy the identity and associativity axioms. For more details about operads, we refer to~\cite{STAN}.
	Let $\mathcal{P}(n)_{\mathbb{S}_n}$ be the space of coinvariants of the operad $\mathcal{P}$, that is, $\mathcal{P}(n)_{\mathbb{S}_n}=\mathcal{P}_n/\{v\sigma-v\mid \sigma\in \mathbb{S}_n,v\in \mathcal{P}_n\}$. Define
	\[G(\mathcal{P})=\prod_{n=1}^{+\infty}\mathcal{P}(n)_{\mathbb{S}_n}.\] It is obvious that we can define the operator $+\colon G(\mathcal{P})\times G(\mathcal{P})\to G(\mathcal{P})$~by
	\[(\overline{a}+\overline{b})_n=\overline{a_n+b_n},\qquad \forall a,b\in G(\mathcal{P}),\quad\forall n\ge 1.\] Indeed, $(G(\mathcal{P}),+)$ is an abelian group.
	We can define another operator $\times\colon G(\mathcal{P})\times G(\mathcal{P})\to G(\mathcal{P})$~by
	\[(\overline{a}\times \overline{b})_n=\sum_{k=1}^{n}\sum_{t_1+t_2+\cdots+t_k=n}\overline{\gamma(b_k;a_{t_1},a_{t_2},\dots,a_{t_k})},\qquad\forall a,b\in \mathcal{P},\quad\forall n\ge 1.\] It follows from \cite[Proposition~4.1]{CPT} that $\times$ is associative and that each left multiplication with respect to $\times$ is an endomorphism of the abelian group $(G(\mathcal P),+)$. Hence $(G(\mathcal P),+,\times,\Phi)$ is a weak twisted post-group, where $\Phi\colon G(\mathcal{P})\to G(\mathcal{P})$ is the trivial map, that is, $\Phi(\overline{a})=(\overline{0},\overline{0},\dots,\overline{0},\cdots)$ for any $a\in \mathcal{P}$.
\end{Example}		
Rings also provide an important class of weak twisted post-groups.
\begin{Example}
	Let $R$ be a ring with multiplication $\cdot$. Then $(R,+,\cdot,\Phi)$ is a weak twisted post-group, where $\Phi$ is given by $\Phi(a)=0$ for any $a\in R$. Moreover, if $R$ is a field, then $(R,+,\cdot,\Phi)$ is a weak twisted post-group but not a twisted post-group as the left multiplication $L^{\rhd}_0$ is not an automorphism.
\end{Example}
Throughout this paper, given a twisted post-group $(G,\cdot,\rhd,\Phi)$, unless otherwise cited,
we will denote by \begin{equation*}e_{a}=(L_{a}^{\rhd})^{-1}\bigl(\Phi(a)^{-1}\cdot a\bigr),
\end{equation*} and
\begin{equation}\label{a+}
	a^{\dag}=(L_{a}^{\rhd})^{-1}\bigl(\Phi(a)^{-1}\cdot e_a\bigr),
\end{equation} for any $a\in G$.

\begin{Definition}
	A homomorphism between (weak) twisted post-groups $(G,\cdot_{G},\rhd_{G},\Phi_{G})$ and $(H,\cdot_{H},\rhd_{H},\Phi_{H})$ is a map $F\colon G\to H$ such that for any $a,b\in G$,
	\begin{subequations}
		\begin{gather}
			F(a)\cdot_{H} F(b)=F(a\cdot_{G} b),\label{TPHM1}\\
			F(a)\rhd_{H}F(b)=F(a\rhd_{G} b),\label{TPHM2}\\
			\Phi_{H}F(a)=F(\Phi_{G}(a)).\label{TPHM3}
		\end{gather}
	\end{subequations}
	
\end{Definition}

\begin{Lemma}\label{LC}
	Let $(G,\cdot,\rhd,\Phi)$ be a twisted post-group with the sub-adjacent operation $\circ$. Then $G$ is right divisible with respect to $\circ$, that is, for any $a,b\in G$, there is a unique $c\in G$, such that $a\circ c=b$.
\end{Lemma}
\begin{proof}
It is not hard to verify by the definition. Indeed, $a\circ c=b$ is equivalent to
$\Phi(a)\cdot(a\rhd c)=b$,
or equivalently
$a\rhd c=\Phi(a)^{-1}\cdot b$.
Since $L_a^\rhd$ is an automorphism of $(G,\cdot)$, the unique solution is
$c=(L_a^\rhd)^{-1}\bigl(\Phi(a)^{-1}\cdot b\bigr)$.
\end{proof}

The next theorem generalises \cite[Theorem~2.4]{BM}.

\begin{Theorem}\label{th1}Let $(G,\cdot_G,\rhd_G,\Phi_G)$ be a twisted post-group with the sub-adjacent operation $\circ_G$. Then we have
	\begin{enumerate}\itemsep=0pt
		\item[$(a)$] $(G,\circ_G)$ is a semigroup. It is called the \textit{sub-adjacent semigroup} of $(G,\cdot_G,\rhd_G,\Phi_G)$;
		\item[$(b)$] For any $a\in G$, $a\circ_G e_a=a$, $a\circ_G a^{\dag}=e_a$, and $\Phi_G(e_a)=1_G$;
		\item[$(c)$] The left multiplication $L^{\rhd_G}$ is a semigroup action of $(G,\circ_G)$ on $(G,\cdot_G)$;
		\item[$(d)$] Let $F\colon (G,\cdot_{G},\rhd_{G},\Phi_{G})\to (H,\cdot_{H},\rhd_{H},\Phi_{H})$ be a twisted post-group homomorphism. Then $F$ induces a semigroup homomorphism from $(G,\circ_{G})$ to $(H,\circ_{H})$.
	\end{enumerate}
\end{Theorem}
\begin{proof}
(a) For any $a,b,c\in G$, we have
\begin{align*}
	(a\circ_G b)\circ_G c =\Phi_{G}(a\circ_G b)\cdot_G((a\circ_G b)\rhd_G c)
 \overset{\eqref{TPG2}}{=}\Phi_G(a\circ_G b)\cdot_G(a\rhd_G (b\rhd_G c)).
\end{align*}
And we have
\begin{align*}
	a\circ_G (b\circ_G c)&=a\circ_G(\Phi_G(b)\cdot_G(b\rhd_G c))
	 =\Phi_G(a)\cdot_G(a\rhd_G (\Phi_G(b)\cdot_G(b\rhd_G c)))\\
\overset{\eqref{TPG1}}&{=}\Phi_G(a)\cdot_G(a\rhd_G\Phi_G(b))\cdot_G(a\rhd_G(b\rhd_G c))\\
	&=(a\circ_G \Phi_G(b))\cdot_G(a\rhd_G(b\rhd_G c))\\
\overset{\eqref{TPG3}}&{=}\Phi_{G}(a\circ_G b)\cdot_G(a\rhd_G(b\rhd_G c))
\end{align*}
This proves the associativity of $\circ_G$.

(b) By the definition of $e_a$,
\[
a\circ e_a
=\Phi(a)\cdot(a\rhd e_a)
=\Phi(a)\cdot \Phi(a)^{-1}\cdot a
=a.
\]
Similarly, by the definition of $a^\dag$,
\[
a\circ a^\dag
=\Phi(a)\cdot(a\rhd a^\dag)
=\Phi(a)\cdot \Phi(a)^{-1}\cdot e_a
=e_a.
\]
Moreover,
\[
a\circ \Phi(e_a)=\Phi(a\circ e_a)=\Phi(a)=a\circ 1_G,
\]
where $a\circ 1_G=\Phi(a)$ since $L_a^\rhd$ is a group automorphism and hence $a\rhd 1_G=1_G$. By Lemma $\ref{LC}$, we obtain $\Phi(e_a)=1_G$.

(c) By \eqref{TPG2}, $L^{\rhd_G}$ defines a semigroup action from $(G,\circ_G)$ to $(G,\cdot_G)$.

(d) For any $a,b\in G$, we have
\begin{align*}
F(a\circ_{G}b)&=F(\Phi_G(a)\cdot_{G}(a\rhd_{G} b))\overset{\eqref{TPHM1}}{=}F(\Phi_G(a))\cdot_{H}F(a\rhd_{G} b)\\
\overset{\eqref{TPHM2},\,\eqref{TPHM3}}&{=}\Phi_H(F(a))\cdot_{H}(F(a)\rhd_{H}F(b))=F(a)\circ_{H} F(b).
\end{align*}
Therefore, $F$ induces a homomorphism of semigroups from $(G,\circ_G)$ to $(H,\circ_H)$.
\end{proof}

\section{Internal structure of twisted post-groups}\label{3}
In this section, we study the internal structure of twisted post-groups.
We generalise the algebraic properties of Rota--Baxter systems of groups
established in \cite[Section 4]{ZH} to twisted post-groups.

Throughout this section, given a twisted post-group $(G,\cdot,\rhd,\Phi)$ with the sub-adjacent operation $\circ$, we denote by $G_a$ the subset
$G_a:=\{b\circ e_a\mid b\in G\}$
of~$G$. As is the case for Rota--Baxter systems of groups, in the next proposition, we show that~$G_a$ is a group with respect to~$\circ$.

\begin{Proposition}\label{SG}
	Let $(G,\cdot,\rhd,\Phi)$ be a twisted post-group and $a\in G$. Then $(G_a,\circ)$ is a group with identity element $e_a$.
\end{Proposition}
\begin{proof}
By Theorem~\ref{th1}, we have \[(b\circ e_a)\circ (c\circ e_a)=(b\circ e_a\circ c)\circ e_a ,\qquad \forall b,c\in G.\] This implies $G_a$ is closed under $\circ$.

For any $a,b\in G$, using \eqref{TPG2} and Theorem~\ref{th1}, we have
\[a\rhd b=(a\circ e_a)\rhd b=a\rhd(e_a\rhd b).\]
Since $L_a^\rhd$ is injective, it follows that $e_a\rhd b=b$. Again using Theorem~\ref{th1}, we have \begin{equation}\label{IDT}
	e_a\circ b=\Phi(e_a)\cdot(e_a\rhd b)=1\cdot b=b.
\end{equation}
Moreover, by \eqref{IDT}, we obtain \[(a\circ e_a)\circ e_a=a\circ(e_a\circ e_a)=a\circ e_a,\] this proves that $e_a$ is the identity element of $G_a$ with respect to $\circ$.

Finally, we show that $G_a$ is closed under the inverse map. For any $b\in G$, we have $b\circ e_a\in G_a$. Using \eqref{IDT} and Theorem~\ref{th1}, we have
\begin{align*}
	(b\circ e_a)\circ	\bigl(b^{\dag}\circ e_a\bigr)&=b\circ \bigl(e_a\circ b^{\dag}\bigr)\circ e_a=b\circ b^{\dag}\circ e_a=\bigl(b\circ b^{\dag}\bigr)\circ e_a=e_b\circ e_a=e_a.
\end{align*}
By Theorem~\ref{th1} and the above identity, we have
\begin{align*}
	b\circ e_a&=(e_a\circ b)\circ e_a=\bigl(b\circ b^{\dag}\circ e_a \circ b\bigr)\circ e_a=b\circ\bigl(b^{\dag}\circ e_a\circ b\circ e_a\bigr),
\end{align*}
it follows from Lemma~\ref{LC} that $e_a= \bigl(b^{\dag}\circ e_a\bigr)\circ(b\circ e_a)$, this proves that $b\circ e_a$ is invertible with respect to $\circ$.
\end{proof}

\begin{Lemma}\label{Phi}
	Let $(G,\cdot,\rhd,\Phi)$ be a twisted post-group with the sub-adjacent operation $\circ$. Then $\Phi(G)=G_1$. Furthermore, for any $a$ in group~$G_1$, its inverse element is~$a^{\dag}$.
\end{Lemma}
\begin{proof}
For any $a\in G$, we have \[\Phi(a)=a\circ 1=a\circ (1\circ e_1)=(a\circ 1)\circ e_1.\] This means $\Phi(G)\subseteq G_{1}$.

Conversely, using Theorem~\ref{th1}, we have \[1\circ e_1=1=e_1\circ 1=\bigl(1\circ 1^{\dag}\bigr)\circ1=1\circ\bigl(1^{\dag}\circ 1\bigr),\] then by Lemma~\ref{LC}, we have $e_1=1^{\dag}\circ 1$. Then we have{\samepage \[a\circ e_1=a\circ \bigl(1^{\dag}\circ 1\bigr)=\bigl(a\circ 1^{\dag}\bigr)\circ 1=\Phi\bigl(a\circ 1^{\dag}\bigr).\] This implies $G_{1}\subseteq \Phi(G)$. We prove that $(\Phi(G),\circ)=(G_1,\circ)$.}

For any $a\in \Phi(G)$, we show that $e_a=e_1$. By the definition we have $a\circ e_a=a$ and by Proposition~\ref{SG}, we have $a\circ e_1=a$. Then using Lemma~\ref{LC}, we have $e_a=e_1$.
It follows from the definition that $a\circ a^{\dag}=e_a=e_1$. Next using Proposition~\ref{SG}, we have
\[a\circ e_1=a=e_1\circ a=e_a\circ a=\bigl(a\circ a^{\dag}\bigr)\circ a=a\circ\bigl(a^{\dag}\circ a\bigr).\] Then by Lemma~\ref{LC}, we have $e_1=a^{\dag}\circ a$, hence the inverse element of $a$ is $a^{\dag}$.
\end{proof}

The group $(\Phi(G),\circ)=(G_1,\circ)$ is called the \textit{sub-adjacent group} of $(G,\cdot,\rhd,\Phi)$.

Immediately from \eqref{TPHM3} and Theorem~\ref{th1}, one deduces that a twisted post-group homomorphism induces a group homomorphism between the sub-adjacent groups.
\begin{Proposition}
	Let $F$ be a twisted post-group homomorphism from $(G,\cdot_G,\rhd_G,\Phi_G)$ to $(H,\cdot_H,\allowbreak \rhd_H,\Phi_H)$. Then the restriction of $F$ on $\Phi_G(G)$ is a group homomorphism from the sub-adjacent group $(\Phi_G(G),\circ_G)$ to the sub-adjacent group $(\Phi_H(H),\circ_H)$.
\end{Proposition}
\begin{proof}
By \eqref{TPHM3}, $F$ maps $\Phi_G(G)$ to $\Phi_H(H)$.
Then using Theorem~\ref{th1}, it follows that $F$ is compatible with $\circ$.
\end{proof}

Now we show that $G_a$ is isomorphic to each other for any $a\in G$ and $G$ is the disjoint union of $G_a$. The next theorem generalises \cite[Lemma 4.2 and Proposition 4.6]{ZH}.
\begin{Theorem}\label{th2}
Let $(G,\cdot,\rhd,\Phi)$ be a twisted post-group with the sub-adjacent operation $\circ$. With the above notations, we have the following statements:
\begin{enumerate}\itemsep=0pt
		\item For any $a,b\in G$, the operator $\pi_{a,b}\colon (G_a,\circ)\to (G_b,\circ)$ given by
		\[\pi_{a,b}(t)=t\circ e_b,\qquad\forall t\in G_a,\] is a group isomorphism.
		\item $G=\bigsqcup G_a$, that is, $G$ is the disjoint union of the groups $G_a$, where $a\in G$.
\end{enumerate}
\end{Theorem}
\begin{proof}
(a) For any $s,t\in G_a$, by Theorem~\ref{th1}, we have
\begin{align*}
	\pi_{a,b}(s)\circ\pi_{a,b}(t) =(s\circ e_b)\circ (t\circ e_b)=s\circ (e_b\circ t)\circ e_b
 =s\circ t\circ e_b=\pi_{a,b}(s\circ t).
\end{align*}
This implies $\pi_{a,b}$ is a group homomorphism.

It remains to show it is bijective.
If there is a $t\in G_a$ such that $\pi_{a,b}(t)=e_b$, it follows from Proposition~\ref{SG} and \eqref{IDT} that
\[t=t\circ e_a=t\circ e_b\circ e_a=\pi_{a,b}(t)\circ e_a=e_b\circ e_a=e_a.\] Hence $\pi_{a,b}$ is injective.
For any $t\in G_b$, again by Proposition~\ref{SG} and \eqref{IDT}, we have \[\pi_{a,b}(t\circ e_a)=(t\circ e_a)\circ e_b=t\circ(e_a\circ e_b)=t\circ e_b=t,\] this proves that $\pi_{a,b}$ is surjective.

(b) For any $a\in G$, by the definition of $e_a$, we have
$a\circ e_a=a$.  This implies $a\in G_a$. If $t\in G_a\cap G_b$ for $a,b\in G$, then we have $t\circ e_a=t=t\circ e_b$ by Proposition~\ref{SG}, it follows from Lemma~\ref{LC} that $e_a=e_b$, which implies $G_a=G_b$. This proves that $G=\bigsqcup G_a$.\end{proof}

Let $(G,\cdot,\rhd,\Phi)$ be a twisted post-group with the sub-adjacent operation $\circ$. Denote the set $\{e_a\mid a\in G\}$ by $K$. It is easy to see that $K$ is a semigroup with respect to $\circ$ by \eqref{IDT}.

Here is the main result of this section. In the next theorem, we generalise the decomposition theorem of Rota--Baxter systems of groups \cite[Theorem 4.7]{ZH} to the case for twisted post-groups. More precisely, we show that as a semigroup, the sub-adjacent semigroup $(G,\circ)$ can be decomposed as the direct product of the sub-adjacent group $(G_1,\circ)$ and the semigroup $(K,\circ)$.

\begin{Theorem}\label{DU}
	Let $(G,\cdot,\rhd,\Phi)$ be a twisted post-group. With the above notations, we have
$(G,\circ)\simeq (G_{1},\circ)\oplus (K,\circ)$ as semigroups.
\end{Theorem}
\begin{proof}
Define $\Psi\colon (G,\circ)\to (G_{1},\circ)\oplus( K,\circ)$ by
$\Psi(a)=(a\circ e_{1},e_a)$, $\forall a\in G$. It is easy to verify that $\Psi$ is well-defined.

First, let us show that $\Psi$ is a semigroup homomorphism.
By the proof of Proposition~\ref{SG}, we have \[\Psi(a\circ b)=(a\circ b\circ e_{1},e_{a\circ b})=(a\circ (e_{1}\circ b)\circ e_{1},e_{a\circ b})=((a\circ e_{1} )\circ (b\circ e_{1}),e_{a\circ b}).\] As $a\circ b\circ e_{a\circ b}=a\circ b$, we have $b\circ e_{a\circ b}=b$ by Lemma~\ref{LC}. And since $b=b\circ e_b$, it follows that $e_{a\circ b}=e_b$. Next, by \eqref{IDT}, we know that $e_b=e_a\circ e_b$. Therefore, we have
\[\Psi(a\circ b)=((a\circ e_{1})\circ (b\circ e_{1}),e_b )=((a\circ e_{1})\circ (b\circ e_{1}),e_a\circ e_b)=\Psi(a)\Psi(b).\] This proves that $\Psi$ is a semigroup homomorphism.

Next, we show $\Psi$ is injective. If there exist $a,b \in G$ such that $a\circ e_{1}=b\circ e_1$ and $e_a=e_b$, then by \eqref{IDT}, we have \[a=a\circ e_a=a\circ (e_{1}\circ e_a)=(a\circ e_{1})\circ e_a=(b\circ e_{1})\circ e_b=b\circ(e_{1}\circ e_b)=b\circ e_b=b,\] this means $\Psi$ is injective.

Finally, we prove $\Psi$ is surjective. For any $a\in G$, since
$a\circ e_b\in G_b$, Theorem~\ref{th2}\,(b) implies $G_{a\circ e_b}=G_b$,
and therefore
$e_{a\circ e_b}=e_b$.

By Theorem~\ref{th2}, for any $a,b\in G$, we have \[\Psi(a\circ e_b)=((a\circ e_b)\circ e_{1},e_{a\circ e_b})=(a\circ(e_b\circ e_{1}),e_b)=(a\circ e_{1},e_b).\] This proves that $\Psi$ is surjective.
\end{proof}

\begin{Remark}
As observed in \cite[Remark~3.3]{JF}, when the base braided monoidal
category is $\mathbf{Set}$, (weak) twisted post-Hopf algebras reduce
to (weak) twisted post-groups. Moreover, idempotents split in
$\mathbf{Set}$, and the image object associated with the idempotent
cocycle $\Phi$ can be identified with
$\operatorname{Im}(\Phi)=\Phi(G)$. Thus, the sub-adjacent group
$(\Phi(G),\circ)$ constructed above may be viewed as the
set-theoretical specialization of the image-object construction
developed in~\cite{JF}.

The more detailed results of this section concerning the groups
$G_a$, their pairwise isomorphisms, the disjoint-union decomposition
of $G$, and the semigroup decomposition in Theorem~\ref{DU} are not
direct consequences of that categorical image construction.
\end{Remark}

\section[(Weak) twisted post-groups, skew trusses, Rota--Baxter systems of groups and braces]{(Weak) twisted post-groups, skew trusses,\\ Rota--Baxter systems of groups and braces}\label{4}

In this section, we study the connections between twisted post-groups and other algebraic systems. The section consists of three parts. In the first part, we show that the category of weak twisted post-groups is isomorphic to the category of skew trusses, and we further consider the skew truss analogue of twisted post-groups. In the second part, we investigate the connections between Rota--Baxter systems of groups and twisted post-groups. Finally, in the third part, we~study the relationship between two-sided twisted post-groups and two-sided braces.

\subsection{(Weak) Twisted post-groups and skew trusses}\label{4.1}

First, let us recall some basic concepts on (skew) braces given in \cite{LG}.

\begin{Definition}
	Let $G$ be a group with two operations $\circ$ and $\cdot$ such that $(G,\cdot)$ and $(G,\circ)$ are both groups and
	\begin{equation}\label{SLB1}
		a\circ(b\cdot c)=(a\circ b)\cdot a^{-1}\cdot(a\circ c),\qquad \forall a,b,c\in G,
	\end{equation}
 where $a^{-1}$ denotes the inverse of $a$ in the group $(G,\cdot)$.
 The triple $(G,\cdot,\circ)$ is called a skew left brace. In particular, if $(G,\cdot)$ is abelian, then $(G,\cdot,\circ)$ is called a left brace.
\end{Definition}
A skew right brace is defined similarly, replacing \eqref{SLB1} by
\[(a\cdot b)\circ c=(a\circ c)\cdot c^{-1}\cdot(b\circ c).\] In particular, if $(G,\cdot)$ is abelian, then $(G,\cdot,\circ)$ is called a right brace. A two-sided (skew) brace is a triple $(G,\cdot,\circ)$ that is both a (skew) left and a (skew) right brace.

Next, let us recall some basic notions of skew trusses given in \cite{TB2}.
\begin{Definition}
	A skew left truss is a quadruple $(G,\cdot,\circ,\Phi)$, where $G$ is a set, $\cdot\colon G\times G\to G$, $\circ\colon G\times G\to G$ and $\Phi\colon G\to G$ are operators such that $(G,\cdot)$ is a group, $(G,\circ)$ is a semigroup and the following identity holds:
	\begin{equation}\label{ST1}
		a\circ (b \cdot c)=(a\circ b)\cdot\Phi(a)^{-1}\cdot(a\circ c),\qquad \forall a, b, c\in G.
	\end{equation}
	The operator $\Phi$ is called the cocycle.
\end{Definition}

A skew right truss is defined in a symmetric way, replacing~\eqref{ST1} by
\begin{equation*}
	(a\cdot b)\circ c=(a\circ c)\cdot\Phi(c)^{-1}\cdot(b\circ c),\qquad \forall a, b, c\in G.
\end{equation*}
$(G,+,\circ,\Phi)$ is called a two-sided truss if $(G,+)$ is an abelian group and it is both a skew left and skew right truss. Without further mention, we will simply call skew left trusses by skew trusses.

Let $(G,\cdot_G,\circ_G,\Phi_G)$ and $(H,\cdot_H,\circ_H,\Phi_H)$ be skew trusses. A skew truss homomorphism between $(G,\cdot_G,\circ_G,\Phi_G)$ and $(H,\cdot_H,\circ_H,\Phi_H)$ is a map $F\colon G\to H$ such that $F$ is both a group homomorphism from $(G,\cdot_G)$ to $(H,\cdot_H)$ and a semigroup homomorphism from $(G,\circ_G)$ to $(H,\circ_H)$.

Now we investigate the relationship between skew trusses and weak twisted post-groups.

\begin{Proposition}\label{EQ1}
Assume $(G,\cdot_G,\rhd_G,\Phi_G)$ is a weak twisted post-group with the sub-adjacent \mbox{operation}~$\circ_G$. Then $(G,\cdot_G,\circ_G,\Phi_G)$ is a skew truss. Moreover, if $F\colon (G,\cdot_G,\rhd_G,\Phi_G)\to \linebreak (H,\cdot_H,\rhd_H,\Phi_H)$ is a weak twisted post-group homomorphism, then $F$ induces a skew truss homomorphism.
\end{Proposition}
\begin{proof}
The associativity of $\circ_G$ follows by the same calculation as in
Theorem~\ref{th1}\,(a), since that calculation only uses that
$L_a^{\rhd_G}$ is an endomorphism of $(G,\cdot_G)$. By Definition \ref{TPG}, we have
\begin{align*}
	(a\circ_G b)\cdot_G \Phi_G(a)^{-1}\cdot_G (a\circ_G c)&=(\Phi_G(a)\cdot_G (a\rhd_G b))\cdot_G \Phi_G(a)^{-1}\cdot_G (\Phi_G(a)\cdot_G (a\rhd_G c))\\&
	=\Phi_G(a)\cdot_G (a\rhd_G b)\cdot_G (a\rhd_G c)\\
\overset{\eqref{TPG1}}&{=}\Phi_G(a)\cdot_G (a\rhd_G (b\cdot_G c))=a\circ_G (b\cdot_G c),
\end{align*} for any $a,b,c\in G$.
This proves that $(G,\cdot_G,\circ_G)$ is a skew truss.
Next, by the definition, we have
\begin{align*}
	F(a\circ_G b)&=F(\Phi_G(a)\cdot_G(a\rhd_G b))=F(\Phi_G(a))\cdot_H F(a\rhd_G b)\\
	&=\Phi_{H}(F(a))\cdot_H (F(a)\rhd_H F(b))=F(a)\circ_H F(b).
\end{align*}
This means $F$ is a skew truss homomorphism.
\end{proof}

\begin{Proposition}\label{EQ2}
	Let $(G,\cdot_G,\circ_G,\Phi_G)$ be a skew truss. Define $\rhd_G\colon G\times G\to G$ by
	\begin{equation*}
		a\rhd_G b=\Phi_G(a)^{-1}\cdot_G(a\circ_G b),\qquad\forall a,b\in G.
	\end{equation*}
	Then $(G,\cdot_G,\rhd_G,\Phi_G)$ is a weak twisted post-group. Moreover, if $F\colon (G,\cdot_G,\circ_G,\Phi_G)\to (H,\cdot_H,\allowbreak \circ_H,\Phi_H)$ is a skew truss homomorphism, then it induces a weak twisted post-group homomorphism.
\end{Proposition}
\begin{proof}
Using \cite[Theorem~2.9]{TB2}, we have \eqref{TPG1} and \eqref{TPG2}. Then by \cite[Lemma~2.3]{TB2}, we obtain~\eqref{TPG3}. Hence $(G,\cdot_G,\rhd_G,\Phi_G)$ is a weak twisted post-group.

Moreover, it follows from \cite[Proposition 2.10]{TB2} that a skew truss homomorphism induces a~weak twisted post-group homomorphism.
\end{proof}

Skew trusses and homomorphisms form a category~$\mathcal{ST}$. Weak twisted post-groups and homomorphisms form a category $\mathcal{WTG}$. The next theorem generalises \cite[Theorem~3.18]{BM}.
\begin{Theorem}\label{th4}
	The two categories $\mathcal{ST}$ and $\mathcal{WTG}$ are isomorphic.
\end{Theorem}	
\begin{proof}
By Proposition~\ref{EQ1}, we can define a functor $\mathcal{F}\colon \mathcal{WTG}\to \mathcal{ST}$, and similarly by Proposition~\ref{EQ2}, we can define a functor $\mathcal{G}\colon \mathcal{ST}\to \mathcal{WTG}$. One can readily verify that $\mathcal{F}\circ \mathcal{G}=\operatorname{Id}$ and $\mathcal{G}\circ \mathcal{F}=\operatorname{Id}$. Therefore, $\mathcal{ST}$ and $\mathcal{WTG}$ are isomorphic.
\end{proof}

\begin{Remark}
Hopf trusses are defined more generally in an arbitrary braided
monoidal category~\cite{RG}. When the underlying category is the category of sets $\mathbf{Set}$, Hopf trusses
are precisely skew trusses \cite[Remark~2.18]{RG}. Hence
Theorem~\ref{th4} can equivalently be viewed as an isomorphism between
the category of weak twisted post-groups and the category of Hopf
trusses in~$\mathbf{Set}$.

By \cite[Theorem~2.20]{RG}, the category of skew trusses is
equivalent to the full subcategory of pointed cosemisimple Hopf
trusses in the category of vector spaces over a field~$\mathbb{F}$.
Consequently, the category of weak twisted post-groups is also
equivalent to this full subcategory.
\end{Remark}

In light of the isomorphism between the category of skew trusses and that of weak twisted post-groups, we now give a~characterization of the skew truss analogue of twisted post-groups.
\begin{Proposition}\label{RD}
	Let $(G,\cdot,\circ,\Phi)$ be a skew truss. Then the following statements are equivalent:
\begin{enumerate}\itemsep=0pt
\item Define $\rhd\colon G\times G\to G$ by
$a\rhd b=\Phi(a)^{-1}\cdot(a\circ b)$, $\forall a,b\in G$,
then $(G,\cdot,\rhd,\Phi)$ is a twisted post-group;
\item The semigroup $(G,\circ)$ is right divisible, that is
		for any $a,b\in G$, there is a unique $t\in G$ such that $a\circ t=b$.
\end{enumerate}
\end{Proposition}

\begin{proof}
The proof follows from \cite[Lemma 3.1]{ZH}.
\end{proof}

\subsection{Twisted post-groups, Rota--Baxter systems of groups and rings}

In this subsection, we study the relationship between Rota--Baxter systems of groups and twisted post-groups.
Let $(G,\cdot)$ be a group and $t\in G$ be an element. Denote the automorphism group of $(G,\cdot)$ by $\Aut(G)$. Recall that a conjugation action by $t$ is an operator $\varphi_t$ which is defined by
$\varphi_t(a)=t^{-1}at$, $a\in G$.  And the set $\{\varphi_t\mid t\in G\}$ under the composition of operators, forms a~subgroup of $\Aut(G)$, called the inner group of $G$. Denote this group by $\Inn(G)$.

On the one hand, every Rota--Baxter system of groups carries the structure of a twisted group. On the other hand, Section~\ref{3} shows that twisted post-groups share similar algebraic properties with Rota--Baxter systems of groups. This raises the following question: under what conditions does a twisted group admit the structure of a Rota--Baxter system of groups? To answer this question, the next proposition provides an equivalent characterization of Rota--Baxter systems of groups under certain conditions.

\begin{Proposition}
Let $(G,\cdot,\rhd,\Phi)$ be a twisted post-group with the sub-adjacent operation $\circ$ such that the center of $G$ is trivial, that is, $Z(G)=\{1\}$. Then the following statements are equivalent:
\begin{enumerate}\itemsep=0pt
\item[$(a)$] There is a Rota--Baxter system of groups $(G,B_1,B_2)$ on $(G,\cdot)$, such that $B_1(a)\cdot B_2(a)=\Phi(a)$ and $L_a^{\rhd}=\varphi_{B_2(a)}$ for any $a\in G$, where $\varphi_a$ is the conjugation action by~$a$;
\item[$(b)$] For any $a\in G$, the left multiplication $L_a^{\rhd}$ lies in $\Inn(G)$.
\end{enumerate}
\end{Proposition}
\begin{proof}
(a)~$\Rightarrow$~(b) It follows from Example~\ref{RG}.

(a)~$\Leftarrow$~(b) As the left multiplication $L_a^{\rhd}$ lies in $\Inn(G)$ for any $a\in G$ and $Z(G)=\{1\}$, there is a unique $t\in G$, such that $L_a^{\rhd}=\varphi_t$. Define the operator $B_2\colon G\to G$ by $B_2(a)=t$ for any $a\in G$. Then define $B_1\colon G\to G$ by
\[
B_1(a)  =\Phi(a)\cdot B_2(a)^{-1},\qquad\forall a\in G.\]
It remains to show that $(G,B_1,B_2)$ is a Rota--Baxter system of groups. By the definition, we have $a\circ b=B_1(a)\cdot b \cdot B_2(a)$. Then using~\eqref{TPG2}, we have
\begin{align*}
& B_2(a)^{-1}\cdot B_2(b)^{-1}\cdot c\cdot B_2(b)\cdot B_2(a)\\
& \qquad =B_2(B_1(a)\cdot b\cdot B_2(a))^{-1}\cdot c\cdot B_2(B_1(a)\cdot b\cdot B_2(a)),\qquad\forall a,b,c\in G.
\end{align*}
 This implies
\begin{equation}\label{RBS}
	B_2(b)\cdot B_2(a)=B_2(B_1(a)\cdot b\cdot B_2(a)).
\end{equation}
Again by \eqref{TPG3}, we have
\begin{align}
	 B_1(B_1(a)\cdot b\cdot B_2(a))\cdot B_2(B_1(a)\cdot b\cdot B_2(a)) &=B_1(a\circ b)\cdot B_2(a\circ b)=\Phi(a\circ b)\notag\\
	&=a\circ \Phi(b)=B_1(a)\cdot B_1(b)\cdot B_2(b)\cdot B_2(a).\!\!\!\label{INT}
\end{align}
It follows from \eqref{RBS} and \eqref{INT} that \[B_1(a)\cdot B_1(b)=B_1(B_1(a)\cdot b\cdot B_2(a)).\] This proves that $(G,B_1,B_2)$ is a Rota--Baxter system of groups.
\end{proof}
\begin{Remark}
For ordinary post-groups, cohomological criteria for determining
whether the post-group structure is induced by a Rota--Baxter
operator have been established in \cite{AC,VG}. A similar
interpretation is suggested by the preceding proposition, although
the twisted case requires some additional care.

Indeed, suppose that $L_a^{\rhd}\in\operatorname{Inn}(G)$ for every
$a\in G$, and choose an element $B_2(a)\in G$ such that
$L_a^{\rhd}=\varphi_{B_2(a)}$.
Since
$L_{a\circ b}^{\rhd}=L_a^{\rhd}\circ L_b^{\rhd}$,
the elements $B_2(a\circ b)$ and $B_2(b)B_2(a)$ induce the same inner
automorphism. Hence
\[
\omega(a,b)
=
B_2(a\circ b)^{-1}B_2(b)B_2(a)
\in Z(G).
\]
The Rota--Baxter system identity for $B_2$ is precisely the condition
$\omega(a,b)=1$ for all $a,b\in G$. Changing the choice of the lifts
$B_2(a)$ by central elements changes $\omega$ by a corresponding
central coboundary. Thus, when $Z(G)$ is nontrivial, the existence of
a Rota--Baxter system should be governed by the vanishing of an
appropriate central cohomological obstruction.

When $Z(G)=\{1\}$, the lifts $B_2(a)$ are unique and this obstruction
vanishes automatically, which explains the centerless assumption in
the preceding proposition. In general, however, $(G,\circ)$ is only
a semigroup, rather than a group. Therefore, a complete analogue of
the cohomological characterization in \cite{AC,VG} would require
a suitable cohomology theory for the associated semigroup or skew
truss. We leave the development of such a theory for future work.
\end{Remark}

In a symmetric way, it is not hard to obtain the following proposition. The proof is similar to the above proposition.
\begin{Proposition}
	Let $(G,\cdot,\lhd,\Phi)$ be a right twisted post-group with the sub-adjacent operation $\circ$ such that the center of $G$ is trivial, that is, $Z(G)=\{1\}$. Then the following statements are equivalent:
\begin{enumerate}\itemsep=0pt
\item[$(a)$] There is a Rota--Baxter system of groups $(G,B_1,B_2)$ on $(G,\cdot)$, such that $B_1(a)B_2(a)=\Phi(a)$ and $R_a^{\lhd}=\varphi_{B_1(a)^{-1}}$ for any $a\in G$, where $\varphi_a$ is the conjugation action by~$a$;
\item[$(b)$] For any $a\in G$, the right multiplication $R_a^{\lhd}$ lies in $\Inn(G)$.
\end{enumerate}
\end{Proposition}

\subsection{Two-sided (weak) twisted post-groups, braces and rings}

In Section~\ref{4.1}, we show that twisted post-groups can be viewed as special cases of skew trusses and characterize their skew truss analogue. Furthermore, \cite[Theorem 5.2]{TB2} states that every two-sided truss corresponds to a nonunital ring. In this subsection, we investigate the ring corresponding to a two-sided twisted post-group.

First, we give the two-sided version of
(weak) twisted post-groups.
\begin{Definition}
	A quintuple $(G,\cdot,\rhd,\lhd,\Phi)$ is called a two-sided (weak) twisted post-group if%
	\begin{enumerate}\itemsep=0pt
		\item $(G,\cdot,\rhd,\Phi)$ is a left (weak) twisted post-group and $(G,\cdot,\lhd,\Phi)$ is a right (weak) twisted post-group;
		\item $(G,\cdot,\rhd,\Phi)$ and $(G,\cdot,\lhd,\Phi)$ have the same sub-adjacent operation $\circ$.
	\end{enumerate}
	
The operator $\Phi$ is called the cocycle of $(G,\cdot,\rhd,\lhd,\Phi)$ and the operator $\circ$ is called the sub-adjacent operation of $(G,\cdot,\rhd,\lhd,\Phi)$. Moreover, if $(G,\cdot)$ is abelian, then $(G,\cdot,\rhd,\lhd,\Phi)$ is called abelian.
\end{Definition}	

\begin{Definition}	A homomorphism of two-sided (weak) twisted post-groups from
$(G,\cdot_{G},\rhd_G,\allowbreak \lhd_G,\Phi_G)$ to $(H,\cdot_{H},\rhd_H,\lhd_H,\Phi_H)$, is a map $F\colon G\to H$, such that $F$ is both a homomorphism of left twisted post-groups from $(G,\cdot_{G},\rhd_G,\Phi_G)$ to $(H,\cdot_{H},\rhd_H,\Phi_H)$ and a homomorphism of right twisted post-groups from $(G,\cdot_{G},\lhd_G,\Phi_G)$ to $(H,\cdot_{H},\lhd_H,\Phi_H)$.
	
\end{Definition}
Denote by $\mathcal{TST}$ the category of two-sided trusses and $\mathcal{AWG}$ the category of abelian two-sided weak twisted post-groups.
Based on Theorem~\ref{th4}, it is not difficult to obtain the following proposition.

\begin{Proposition}\label{TST}
The categories $\mathcal{TST}$ and $\mathcal{AWG}$ are isomorphic.
\end{Proposition}
Next based on Proposition~\ref{SG} and Theorem~\ref{DU}, we discuss the relationship between two-sided (weak) twisted post-groups and nonunital rings.
\begin{Proposition}
Let $(G,+,\rhd,\lhd,\Phi)$ be an abelian two-sided weak twisted post-group with the sub-adjacent operation $\circ$. Then $(G,+,\circ)$ is a nonunital ring if and only if $\Phi$ is trivial, that is, $\Phi\colon G\to G$ is given by $\Phi(a)=0$ for any $a\in G$. Furthermore, let $(G,+,\rhd,\lhd,\Phi)$ be an abelian two-sided twisted post-group with the sub-adjacent operation $\circ$. Then its cocycle $\Phi$ is trivial if and only if $(G,+,\circ)$ is the zero ring, that is, a ring that consists of one element.
\end{Proposition}
\begin{proof}
If $(G,+,\circ)$ is a nonunital ring, then for any $a,b,c\in G$, we have
\begin{align*}
2\Phi(a)+(a\rhd(b+c))& =\Phi(a)+(a\rhd b)+\Phi(a)+(a\rhd c)\\
& =(a\circ b)+(a\circ c)=a\circ(b+c)=\Phi(a)+(a\rhd(b+c)),
\end{align*}
and thus $\Phi$ is trivial.
Conversely, if $(G,+,\rhd,\lhd,\Phi)$ is a two-sided weak twisted post-group such that $\Phi$ is trivial, then we have $\circ=\rhd=\lhd$, and one can easily verify that $(G,+,\circ)$ is a nonunital ring.

If $(G,+,\rhd,\lhd,\Phi)$ is an abelian two-sided twisted post-group such that $\Phi$ is trivial, then $\circ=\rhd=\lhd$ and $(G,+,\circ)$ is a nonunital ring. Then for any $a\in G$, we have $L^{\rhd}_{0}(a)=0\rhd a=0\circ a=0$, and thus $a=0$. This implies that $(G,+,\circ)$ is the zero ring. The converse implication is obvious.\end{proof}

In the next proposition, we show that, for any twisted post-group, there exists a transformation yielding a twisted post-group with an idempotent cocycle.
\begin{Proposition}\label{IDP}
	Let $(G,\cdot,\rhd,\Phi)$ be a twisted post-group with the sub-adjacent operation $\circ$. Define $\blacktriangleright\colon G\times G\to G$ by
$a\blacktriangleright b=\bigl(a\circ 1^{\dag}\bigr)\rhd b$, $\forall a,b\in G$,  and define $\Psi\colon G\to G$ by
$\Psi(a)=\Phi\bigl(a\circ 1^{\dag}\bigr)$. Then $\Psi$ is idempotent and $(G,\cdot,\blacktriangleright,\Psi)$ is a twisted post-group.
\end{Proposition}

\begin{proof}By \eqref{TPG2}, we have
$a\blacktriangleright b=a\rhd \bigl(1^{\dag}\rhd b\bigr)$, $\forall a,b\in G$.  Therefore, $L_{a}^{\blacktriangleright}\colon G\to G$ given by
$L_{a}^{\blacktriangleright}(b)=a\blacktriangleright b$, $\forall b\in G$ is an automorphism of $(G,\cdot)$. Hence \eqref{TPG1} holds.

Define $\bullet\colon G\times G\to G $ by
$a\bullet b=\Psi(a)\cdot(a\blacktriangleright b)$, $\forall a,b\in G$.
Then for any $a,b,c\in G$, we have
\begin{align*}
	(a\bullet b)\blacktriangleright c&=(\Psi(a)\cdot(a\blacktriangleright b))\blacktriangleright c=\bigl(\Phi\bigl(a\circ 1^{\dag}\bigr)\cdot\bigl(\bigl(a\circ 1^{\dag}\bigr)\rhd b\bigr)\bigr)\blacktriangleright c\\
	&=\bigl(\bigl(a\circ 1^{\dag}\bigr)\circ b\bigr)\blacktriangleright c=\bigl(a\circ 1^{\dag}\circ b\circ 1^{\dag}\bigr)\rhd c\\
	&=\bigl(\bigl(a\circ 1^{\dag}\bigr)\circ \bigl(b\circ 1^{\dag}\bigr)\bigr)\rhd c
	\overset{\eqref{TPG2}}{=}\bigl(\bigl(a\circ 1^{\dag}\bigr)\rhd \bigl(b\circ 1^{\dag}\bigr) \rhd c\bigr)\\
	&=a\blacktriangleright (b\blacktriangleright c).
\end{align*}
This proves that \eqref{TPG2} holds.

It remains to show \eqref{TPG3} holds. For any $a,b\in G$,
we have
\begin{align*}
	a\bullet \Psi(b)&=\Psi(a)\cdot(a\blacktriangleright \Psi(b))=\Phi\bigl(a\circ 1^{\dag}\bigr)\cdot\bigl(\bigl(a\circ 1^{\dag}\bigr)\rhd \bigl(\Phi\bigl(b\circ 1^{\dag}\bigr)\bigr)\bigr)\\
	&=\bigl(a\circ 1^{\dag}\bigr)\circ \Phi\bigl(b\circ 1^{\dag}\bigr)\overset{\eqref{TPG3}}{=}\Phi\bigl(\bigl(a\circ 1^{\dag}\bigr)\circ\bigl(b\circ 1^{\dag}\bigr)\bigr)\\
	&=\Phi\bigl(\bigl(\bigl(a\circ 1^{\dag}\bigr)\circ b\bigr)\circ 1^{\dag}\bigr)
	=\Phi\bigl((a\bullet b) \circ 1^{\dag}\bigr)=\Psi(a\bullet b).
\end{align*}	
This implies \eqref{TPG3} holds.

By the definition and Lemma~\ref{Phi}, we have
\begin{align*}
\Psi^{2}(a)& =\Phi\bigl(\Phi\bigl(a\circ 1^{\dag}\bigr)\circ 1^{\dag}\bigr)=a\circ 1^{\dag}\circ 1\circ 1^{\dag}\circ 1\\
& =a\circ e_1\circ e_1=a\circ e_1=a\circ 1^{\dag}\circ 1=\Phi(a\circ 1^{\dag})=\Psi(a).
\end{align*}
This proves that $\Psi$ is idempotent.
\end{proof}
%\begin{Proposition}

%\end{Proposition}	
\begin{Proposition}\label{SKB}
Let $(G,\cdot,\rhd,\Phi)$ be a twisted post-group with
sub-adjacent operation $\circ$. Let~$\blacktriangleright$ and $\Psi$
be the operations defined in Proposition~{\rm \ref{IDP}}, namely,
\[
a\blacktriangleright b=\bigl(a\circ1^{\dagger}\bigr)\rhd b,
\qquad
\Psi(a)=\Phi\bigl(a\circ1^{\dagger}\bigr).
\]
If $\Phi$ is surjective, then $\Psi=\operatorname{id}_G$.
Consequently, $(G,\cdot,\blacktriangleright)$ is a post-group.
Moreover, its sub-adjacent operation $\bullet$ is given by
\[
a\bullet b
=a\cdot(a\blacktriangleright b)
=(a\circ1^{\dagger})\circ b,
\qquad \forall a,b\in G,
\]
and $(G,\cdot,\bullet)$ is a skew left brace.
\end{Proposition}

\begin{proof}
By Proposition~\ref{IDP},
$(G,\cdot,\blacktriangleright,\Psi)$ is a twisted post-group and
$\Psi$ is idempotent. By Theorem~\ref{th1}, we have
$1\circ1^{\dagger}=e_1$,
and by \eqref{IDT},
$e_1\circ1=1$.
Therefore, for every $a\in G$,
\begin{align*}
\Psi(a\circ1)
&=\Phi\bigl((a\circ1)\circ1^{\dagger}\bigr)
=\Phi\bigl(a\circ\bigl(1\circ1^{\dagger}\bigr)\bigr)
=\Phi(a\circ e_1)\\
&=(a\circ e_1)\circ1
=a\circ(e_1\circ1)
=a\circ1
=\Phi(a).
\end{align*}
It follows that
$\operatorname{Im}(\Phi)\subseteq\operatorname{Im}(\Psi)$.
Since $\Phi$ is surjective, $\Psi$ is also surjective. Since $\Psi$
is both surjective and idempotent, it follows that
$\Psi=\operatorname{id}_G$.
Indeed, for any $a\in G$, choose $x\in G$ such that $\Psi(x)=a$.
Then
$\Psi(a)=\Psi^2(x)=\Psi(x)=a$.
Thus $(G,\cdot,\blacktriangleright)$ is a post-group.

Its sub-adjacent operation is
$
a\bullet b
=a\cdot(a\blacktriangleright b)$.
Moreover, since $\Psi(a)=a$, we have
\begin{align*}
a\bullet b
 =\Psi(a)\cdot(a\blacktriangleright b)
 =\Phi\bigl(a\circ1^{\dagger}\bigr)
 \cdot\bigl(\bigl(a\circ1^{\dagger}\bigr)\rhd b\bigr)
 =\bigl(a\circ1^{\dagger}\bigr)\circ b.
\end{align*}
Finally, it follows from \cite[Proposition~3.22]{BM} that
$(G,\cdot,\bullet)$ is a skew left brace.
\end{proof}

Recall from \cite{LG} that a pair $(X,r)$ is called a set-theoretical solution of the Yang--Baxter equation, if $X$ is a set and \[r\colon \ X\times X\to X\times X,\qquad r(a,b)=(\varphi_a(b),\psi_b(a)),\qquad \forall a,b\in X\] is a bijective map, such that
\[(r\times \id)(\id\times r)(r\times \id)=(\id\times r)(r\times \id)(\id\times r).\]
The solution $(X,r)$ is called non-degenerate if $\varphi_a$ and $\psi_b$ are bijective for all $a,b\in  X$.

Using the Proposition~\ref{SKB} and \cite[Theorem 3.1]{LG}, we derive the following corollary.
\begin{Corollary}
	Let $(G,\cdot,\rhd,\Phi)$ be a twisted post-group such that $\Phi$ is surjective. Let $\blacktriangleright\colon G\times G\to G$ and $\bullet\colon G\times G\to G$ be the operators given in Propositions~{\rm \ref{IDP}} and~{\rm \ref{SKB}}. Define $r\colon G\times G\to G\times G$ by
\[r(a,b)=\bigl(L^{\blacktriangleright}_a(b),\bigl(L^{\blacktriangleright}_{L^{\blacktriangleright}_a(b)}\bigr)^{-1}\bigl((a\bullet b)^{-1}\cdot a\cdot(a\bullet b)\bigr)\bigr),\qquad\forall a,b\in G.\] Then $(G,r)$ is a non-degenerate set-theoretical solution of the Yang--Baxter equation. Where for any $a\in G$, the left multiplication $L^{\blacktriangleright}_a$ is given by $L^{\blacktriangleright}_a(b)=a\blacktriangleright b$ for any $b\in G$.
\end{Corollary}

In the next theorem, we show that a two-sided post-group has the structure of two-sided skew brace.
\begin{Theorem}\label{th3}
	Let $(G,\cdot,\rhd,\lhd,\Phi)$ be a two-sided twisted post-group with the sub-adjacent operation $\circ$. Define $\bullet\colon G\times G\to G$ by
$a\bullet b=(a\circ 1^{\dag})\circ b$, $\forall a,b\in G$,  where $1^{\dag}$ is defined in \eqref{a+}. Then $(G,\cdot,\bullet)$ is a two-sided skew brace.
\end{Theorem}
\begin{proof}First, we show that $(G,\circ)$ is a group. By \eqref{IDT}, we have $e_1\circ a=a$ for any $a$. As $(G,\cdot,\lhd,\Phi)$ is a right twisted post-group, in a symmetric way, there is $\bar{e}_1$ such that $a\circ \bar{e}_1=a$. Then $1\circ e_1=1=1\circ \bar{e}_1$ by the definition. It follows from Lemma~\ref{LC} that $e_1=\bar{e}_1$. This implies $e_1$ is the identity element of $(G,\circ)$. Next using Lemma~\ref{Phi}, we have $(G,\circ)$ is a group and $\Phi$ is surjective. Finally, it follows from Proposition~\ref{SKB} that $(G,\cdot,\bullet)$ is a skew left brace. It is similar to prove that $(G,\cdot,\bullet)$ is a skew right brace. Hence $(G,\cdot,\bullet)$ is a two-sided skew brace.
\end{proof}

Conversely, it is easy to see that every two-sided brace has the structure of two-sided twisted post-group.

Recall that a triple $(R,+,\cdot)$ is called a radical ring if $(R,+,\cdot)$ is a nonunital ring such that for any $a\in R$, there is a $b\in R$ such that $a+b+a\cdot b=0$. In the next corollary, we show that every two-sided twisted post-group is associated with a radical ring.

\begin{Corollary}\label{ring}
Let $(G,+,\rhd,\lhd,\Phi)$ be a two-sided twisted post-group with the sub-adjacent operation $\circ$ such that $(G,+)$ is abelian. Let $\bullet$ be the operator given in Theorem~{\rm \ref{th3}}. Define $*\colon G\times G\to G$ by
$a*b=a\bullet b-a-b$, $\forall a,b\in G$.  Then $(G,+,*)$ is a radical ring.
\end{Corollary}
\begin{proof}
The proof follows from Theorem~\ref{th3} and \cite[Proposition 1]{Cedo}.
\end{proof}

\section{Differentiation of twisted post-Lie groups}\label{5}
In this section, we first introduce the notion of twisted post-Lie algebras and show that each twisted post-Lie algebra gives rise to another Lie algebra, called the sub-adjacent Lie algebra. Next, we define twisted post-Lie groups and study their differentiations. In particular, we prove that the differentiation of any twisted post-Lie group yields a twisted post-Lie algebra. Finally, we show that the differentiation of the sub-adjacent Lie group of a twisted post-Lie group coincides with the sub-adjacent Lie algebra of its corresponding twisted post-Lie algebra.

\begin{Definition}
	A twisted post-Lie algebra $(\mathfrak{g},[\cdot,\cdot],\triangleright,\varphi)$ consists of a Lie algebra $(\mathfrak{g},[\cdot,\cdot])$ with linear operators $\triangleright\colon \mathfrak{g}\otimes \mathfrak{g}\to \mathfrak{g}$ and $\varphi\colon \mathfrak{g}\to \mathfrak{g}$ such that for any $x,y,z\in \mathfrak{g}$,
\begin{subequations}
		\begin{gather}
			x\triangleright [y,z]=[x\triangleright y,z]+[y,x\triangleright z];\label{TPL1}\\
			([\varphi(x),\varphi(y)]+x\triangleright \varphi(y)-y\triangleright \varphi(x))\triangleright z=x\triangleright(y\triangleright z)-y\triangleright(x\triangleright z);\label{TPL2}\\
			\varphi\left(\frac{1}{2}[\varphi(x),y]+\frac{1}{2}[x,\varphi(y)]+x\triangleright y-y\triangleright x\right)=[\varphi(x),\varphi(y)]+x\triangleright \varphi(y)-y\triangleright \varphi(x);\label{TPL3}\\	
			\varphi(x)\triangleright y=x\triangleright y;\label{TPL4}\\
			\varphi^2(x)=\varphi(x).\label{TPL5}
		\end{gather}
	\end{subequations}
	
	The operator $[\cdot,\cdot]_{\triangleright}\colon \mathfrak{g}\times \mathfrak{g}\to \mathfrak{g}$ defined by
	\[[x,y]_{\triangleright}=[\varphi(x),\varphi(y)]+x\triangleright \varphi(y)-y\triangleright \varphi(x),\qquad\forall x,y\in\mathfrak{g},\] is called the sub-adjacent operation of $\mathfrak{g}$.
\end{Definition}

Note that if $\varphi$ is the identity map, then $(\mathfrak{g},[\cdot,\cdot],\triangleright)$ is a post-Lie algebra. First we show that~$\mathfrak{g}$ forms a Lie algebra with respect to the sub-adjacent operation $[\cdot,\cdot]_{\triangleright}$.
\begin{Proposition}
Let $(\mathfrak{g},[\cdot,\cdot],\triangleright,\varphi)$ be a twisted post-Lie algebra with the sub-adjacent operation~$[\cdot,\cdot]_{\triangleright}$. Then $\mathfrak{g}$ is a Lie algebra with respect to $[\cdot,\cdot]_{\triangleright}$.
\end{Proposition}

\begin{proof}
It follows from~\eqref{TPL4} that
\begin{gather}\label{Lie}
	[x,y]_{\triangleright}=[\varphi(x),\varphi(y)]+x\triangleright \varphi(y)-y\triangleright\varphi(x)=[\varphi(x),\varphi(y)]+\varphi(x)\triangleright \varphi(y)-\varphi(y)\triangleright\varphi(x).
\end{gather}
And it follows from \eqref{TPL3} and \eqref{TPL5} that
\begin{align*}
	&[\varphi(x),\varphi(y)]+\varphi(x)\triangleright \varphi(y)-\varphi(y)\triangleright\varphi(x)\\
&\qquad =\bigl[\varphi^2(x),\varphi^2(y)\bigr]+\varphi(x)\triangleright \varphi^2(y)-\varphi(y)\triangleright\varphi^2(x)\\
&\qquad	=\varphi\left(\frac{1}{2}\bigl[\varphi^2(x),\varphi(y)\bigr]+\frac{1}{2}\bigl[\varphi(x),\varphi^2(y)\bigr]+\varphi(x)\triangleright\varphi(y) -\varphi(y)\triangleright\varphi(x)\right).
\end{align*}
Again using \eqref{TPL5}, we have
\begin{align*}
	&[\varphi(x),\varphi(y)]+\varphi(x)\triangleright \varphi(y)-\varphi(y)\triangleright\varphi(x)\\
&\qquad	=\varphi\left(\frac{1}{2}[\varphi(x),\varphi(y)]+\frac{1}{2}[\varphi(x),\varphi(y)]+\varphi(x)\triangleright\varphi(y) -\varphi(y)\triangleright\varphi(x)\right)\\
&\qquad	=\varphi([\varphi(x),\varphi(y)]+\varphi(x)\triangleright\varphi(y)-\varphi(y)\triangleright\varphi(x))
\end{align*} for any $x,y\in\mathfrak{g}$. This implies{\samepage
\begin{equation}\label{Lie1}
	[x,y]_{\triangleright}=\varphi([x,y]_{\triangleright})
\end{equation} by \eqref{Lie}.}

Then for any $x,y,z\in \mathfrak{g}$, by \eqref{TPL2}, \eqref{TPL4}, \eqref{Lie} and \eqref{Lie1}, we have
\begin{align*}
&[x,[y,z]_{\triangleright}]_{\triangleright}+[y,[z,x]_{\triangleright}]_{\triangleright} +[z,[x,y]_{\triangleright}]_{\triangleright}\\
&	\quad=\varphi(x)\triangleright \varphi([y,z]_{\triangleright})+\varphi(y)\triangleright \varphi([z,x]_{\triangleright})+\varphi(z)\triangleright \varphi([x,y]_{\triangleright})-\varphi([y,z]_{\triangleright})\triangleright \varphi(x)\\
&\quad\quad{} -\varphi([z,x]_{\triangleright})\triangleright \varphi(y)-\varphi([x,y]_{\triangleright})\triangleright \varphi(z)
+[\varphi(x),\varphi([y,z]_{\triangleright})]+[\varphi(y),\varphi([z,x]_{\triangleright})]\\
&\quad\quad{} +[\varphi(z),\varphi([x,y]_{\triangleright})]\\
	&\quad=\varphi(x)\triangleright [y,z]_{\triangleright}+\varphi(y)\triangleright [z,x]_{\triangleright}+\varphi(z)\triangleright [x,y]_{\triangleright}-[y,z]_{\triangleright}\triangleright \varphi(x)\\
&\quad\quad{}-[z,x]_{\triangleright}\triangleright \varphi(y)-[x,y]_{\triangleright}\triangleright \varphi(z)
	+[\varphi(x),[y,z]_{\triangleright}]+[\varphi(y),[z,x]_{\triangleright}]+[\varphi(z),[x,y]_{\triangleright}]\\
&	\quad=\varphi(x)\triangleright (\varphi(y)\triangleright \varphi(z))-\varphi(x)\triangleright(\varphi(z)\triangleright \varphi(y))+\varphi(y)\triangleright(\varphi(z)\triangleright \varphi(x))\\
&\quad\quad{}-\varphi(y)\triangleright(\varphi(x)\triangleright \varphi(z))+\varphi(z)\triangleright (\varphi(x)\triangleright \varphi(y))-\varphi(z)\triangleright(\varphi(y)\triangleright \varphi(x))\\
&\quad\quad {}+\varphi(x)\triangleright [\varphi(y),\varphi(z)]+\varphi(y)\triangleright [\varphi(z),\varphi(x)]+\varphi(z)\triangleright [\varphi(x),\varphi(y)]
	\\
&\quad\quad{} -\varphi(y)\triangleright(\varphi(z)\triangleright \varphi(x))+\varphi(z)\triangleright(\varphi(y)\triangleright \varphi(x))-\varphi(z)\triangleright (\varphi(x)\triangleright \varphi(y))\\
&\quad\quad{} +\varphi(x)\triangleright(\varphi(z)\triangleright \varphi(y))-\varphi(x)\triangleright(\varphi	(y)\triangleright \varphi(z))+\varphi(y)\triangleright(\varphi(x)\triangleright \varphi(z))\\
&\quad\quad{} +[\varphi(x),[y,z]_{\triangleright}]+[\varphi(y),[z,x]_{\triangleright}]+[\varphi(z),[x,y]_{\triangleright}].
\end{align*}
Next, by \eqref{TPL1}, \eqref{Lie1} and the Jacobi identity, we have
\begin{align*}
&[x,[y,z]_{\triangleright}]_{\triangleright}+[y,[z,x]_{\triangleright}]_{\triangleright} +[z,[x,y]_{\triangleright}]_{\triangleright}\\
&\quad=\varphi(x)\triangleright [\varphi(y),\varphi(z)]+\varphi(y)\triangleright [\varphi(z),\varphi(x)]+\varphi(z)\triangleright [\varphi(x),\varphi(y)]+[\varphi(x),[y,z]_{\triangleright}]\\
&\quad\quad{}+[\varphi(y),[z,x]_{\triangleright}]+[\varphi(z),[x,y]_{\triangleright}]\\
&	\quad =[\varphi(x)\triangleright \varphi(y),\varphi(z)]+[\varphi(y),\varphi(x)\triangleright \varphi(z)]+[\varphi(y)\triangleright \varphi(z),\varphi(x)]+[\varphi(z),\varphi(y)\triangleright \varphi(x)]\\
&\quad\quad{} +[\varphi(z)\triangleright \varphi(x),\varphi(y)]+[\varphi(x),\varphi(z)\triangleright \varphi(y)]+
	[\varphi(x),\varphi(y)\triangleright \varphi(z)]-[\varphi(x),\varphi(z)\triangleright \varphi(y)]\\
&\quad\quad{}+[\varphi(x),[\varphi(y),\varphi(z)]]+[\varphi(y),\varphi(z)\triangleright \varphi(x)]-[\varphi(y),\varphi(x)\triangleright \varphi(z)]\\
&\quad\quad{}+[\varphi(y),[\varphi(z),\varphi(x)]]+[\varphi(z),\varphi(x)\triangleright \varphi(y)]-[\varphi(z),\varphi(y)\triangleright \varphi(x)]+[\varphi(z),[\varphi(x),\varphi(y)]]\\
&	\quad=[\varphi(x),[\varphi(y),\varphi(z)]]+[\varphi(y),[\varphi(z),\varphi(x)]]+[\varphi(z),[\varphi(x),\varphi(y)]]=0.\tag*{\qed}
\end{align*}\renewcommand{\qed}{}
\end{proof}

Then we show that $\varphi(\mathfrak{g})$ is a Lie algebra with respect to $[\cdot,\cdot]_{\triangleright}$.
\begin{Proposition}
	Let $(\mathfrak{g},[\cdot,\cdot],\triangleright,\varphi)$ be a twisted post-Lie algebra. Let $(\mathfrak{g},[\cdot,\cdot]_{\triangleright})$ be its sub-adjacent Lie algebra. Then $\varphi(\mathfrak{g})\subseteq \mathfrak{g}$ is a Lie subalgebra with respect to $[\cdot,\cdot]_{\triangleright}$.
\end{Proposition}
\begin{proof}
Using \eqref{TPL3} and \eqref{TPL4}, for any $x,y\in\mathfrak{g}$, we have
\begin{align*}
	[\varphi(x),\varphi(y)]_{\triangleright}
	&=\bigl[\varphi^{2}(x),\varphi^{2}(y)\bigr]+\varphi(x)\triangleright \varphi^{2}(y)-\varphi(y)\triangleright\varphi^{2}(x)\\
	&=\varphi\left(\frac{1}{2}\bigl[\varphi^{2}(x),\varphi(y)\bigr]+\frac{1}{2}\bigl[\varphi(x),\varphi^{2}(y)\bigr] +\varphi(x)\triangleright\varphi(y)-\varphi(y)\triangleright \varphi(x)\right).
\end{align*}
This proves the assertion.
\end{proof}

Now we give the definition of homomorphisms between twisted post-Lie algebras.
\begin{Definition}A homomorphism from a twisted post-Lie algebra $(\mathfrak{g},[\cdot,\cdot]_{\mathfrak{g}},\triangleright_{\mathfrak{g}},\varphi_{\mathfrak{g}})$ to $(\mathfrak{h},[\cdot,\cdot]_{\mathfrak{h}},\allowbreak \triangleright_{\mathfrak{h}},\varphi_{\mathfrak{h}})$ is a linear map $f$ satisfying the following equalities:
	\begin{equation*}
		f([x,y]_{\mathfrak{g}})=[f(x),f(y)]_{\mathfrak{h}},
\qquad
		f(x\triangleright_{\mathfrak{g}}y)=f(x)\triangleright_{\mathfrak{h}}f(y),
\qquad
		f(\varphi_{\mathfrak{g}}(x))=\varphi_{\mathfrak{h}}(f(x)),
	\end{equation*}
	for any $x,y\in\mathfrak{g}$.
\end{Definition}

Next, we introduce two lemmas on twisted post-groups, which will be used later.
\begin{Lemma}\label{IDM1}
	Let $(G,\cdot,\rhd,\Phi)$ be a twisted post-group. Then $\Phi$ is idempotent if and only if $\Phi(1)=1$.
\end{Lemma}
\begin{proof}
If $\Phi$ is idempotent, then we have $1\circ 1\circ 1=\Phi^{2}(1)=\Phi(1)=1\circ 1$. Then, by Lemma~\ref{LC}, we have $\Phi(1)=1\circ 1=1$.

Conversely, if $\Phi(1)=1$, then for any $a\in G$, we have
\[\Phi^{2}(a)=a\circ 1\circ 1=a\circ \Phi(1)=a\circ 1=\Phi(a),\] which implies $\Phi$ is idempotent.
\end{proof}

\begin{Lemma}\label{IDM2}
	Let $(G,\cdot,\rhd,\Phi)$ be a twisted post-group with the sub-adjacent operation $\circ$ such that~$\Phi$ is idempotent. Then the following equality holds:
	\begin{equation}\label{FRT}
		a\rhd b=\Phi(a)\rhd b,\qquad\forall a,b\in G.
	\end{equation}
\end{Lemma}
\begin{proof}
Using Lemma~\ref{IDM1}, we have $1\circ 1=1=1\circ e_1$. This implies $e_1=1$ by Lemma~\ref{LC}. Then by \eqref{TPG2}, we obtain
\[\Phi(a)\rhd b=(a\circ 1)\rhd b=a\rhd (1\rhd b)=a\rhd (1\circ b)=a\rhd b,\qquad\forall a,b\in G.\tag*{\qed}
\]
\renewcommand{\qed}{}
\end{proof}

Let $(\mathfrak{g},[\cdot,\cdot])$ be the corresponding Lie algebra of the Lie group $(G,\cdot)$. Denote by $\Aut(\mathfrak{g})$ and $\Der(\mathfrak{g})$ the Lie group of automorphisms and the Lie algebra of derivations on the Lie algebra $(\mathfrak{g},[\cdot,\cdot])$. Let
$\exp\colon \mathfrak{g}\to G$,
be its exponential map. The Lie bracket $[\cdot,\cdot]$ and group multiplication are connected by the formula
\[[x,y]=\left.\frac{\d}{\d t} \right|_{t=0}\left.\frac{\d}{\d s} \right|_{s=0}\exp(sx)\exp(ty)\exp(-sx),\qquad \forall x,y\in \mathfrak{g}.\]
Given a map $F\colon G \to G$, denote the tangent map of $F$ by $F_{*1}$, that is,
\[F_{*1}(x)=\left.\frac{\mathrm{d}}{\mathrm{d}t}\right|_{t=0}(F(\exp(tx))),\qquad\forall x\in \mathfrak{g}.\]
As shown in \cite{BM}, since $L^{\rhd}_{a}\in \Aut(G)$, it follows that $(L^{\rhd}_{a})_{*1}\in \Aut(\mathfrak{g})$. Define $\triangleright\colon \mathfrak{g}\otimes\mathfrak{g}\to \mathfrak{g}$ by
\[x\triangleright y=L^{\rhd}_{*1}(x)(y)=\left.\frac{\mathrm{d}}{\mathrm{d}t}\right|_{t=0}\bigl(L^{\rhd}_{\exp(tx)}\bigr)_{\ast 1}(y)=\left.\frac{\d}{\d t} \right|_{t=0}\left.\frac{\d}{\d s} \right|_{s=0}L^{\rhd}_{\exp(tx)}(\exp(sy)).\]
\begin{Definition}\label{def TPG}
A twisted post-Lie group is a twisted post-group
$(G,\cdot,\rhd,\Phi)$ satisfying the following conditions:
\begin{enumerate}\itemsep=0pt
 \item $(G,\cdot)$ is a finite-dimensional real Lie group, and
 $\rhd\colon G\times G\to G$ and $\Phi\colon G\to G$ are smooth;

 \item $\Phi(1)=1$ and $\Phi$ is idempotent, that is,
$ \Phi^2=\Phi$;
 \item Writing
$ \varphi=\Phi_{\ast 1}\colon \mathfrak g\to\mathfrak g$,
 there exists a neighbourhood $U$ of $0$ in $\mathfrak g$ such
 that
$\Phi(\exp x)=\exp(\varphi(x))$, $x\in U$.
\end{enumerate}
\end{Definition}

\begin{Remark} Condition~(c) allows us to avoid taking second derivatives of the $G$-valued map $\Phi\circ\exp$. Indeed, after shrinking $U$ if necessary and denoting by $\log$ the local inverse of $\exp$, condition~(c) gives $(\log\circ\Phi\circ\exp)(x)=\varphi(x)$, $x\in U$.  Thus, all expansions below are carried out in the vector space $\mathfrak g$. In particular, in exponential coordinates,
\[
\begin{aligned}
\left.\frac{\d}{\d t}\right|_{t=0}
\left.\frac{\d}{\d s}\right|_{s=0}
\left(
\log\Phi\bigl(\exp(tx)\circ\exp(sy)\bigr)
\right)
&=
\varphi\left(
\frac{1}{2}[\varphi(x),y]+x\triangleright y
\right),\\
\left.\frac{\d}{\d t}\right|_{t=0}
\left.\frac{\d}{\d s}\right|_{s=0}
\left(
\log\Phi\bigl(\exp(sy)\circ\exp(tx)\bigr)
\right)
&=
\varphi\left(
\frac{1}{2}[\varphi(y),x]+y\triangleright x
\right).
\end{aligned}
\]

 Without condition~(c), the map $\log\circ\Phi\circ\exp$ may have nonlinear second-order terms, which would contribute additional terms to these identities.

Condition (c)
is strictly weaker than requiring $\Phi$ to be a Lie group
homomorphism. Indeed, every Lie group homomorphism satisfies
this identity for all $X\in\mathfrak{g}$, whereas the converse
fails even for smooth idempotent maps. This condition requires
$\Phi$ to be linear in exponential coordinates near the identity,
but does not force $\varphi$ to preserve the Lie bracket.
 \end{Remark}

The following examples illustrate that the cocycle of a
twisted post-Lie group need not be a~Lie group homomorphism,
even when the exponential compatibility holds globally.
\begin{Example}\label{ex:heisenberg-twisted-post-Lie}
Let $G=\mathbb{R}^3$ be the three-dimensional Heisenberg
Lie group with multiplication
\[
(x,y,z)\cdot(u,v,w)
=
\left(
x+u,
y+v,
z+w+\frac12(xv-yu)
\right).
\]
Its identity element is $1=(0,0,0)$, and
$(x,y,z)^{-1}=(-x,-y,-z)$.
Define $\Phi\colon G\to G$ and $\rhd\colon G\times G\to G$ by
$\Phi(x,y,z)=(x,y,0)$ and
\[
(x,y,z)\rhd(u,v,w)
=
\left(
u,
v,
w-\frac12(xv-yu)
\right).
\]
Then $(G,\cdot,\rhd,\Phi)$ is a twisted post-Lie group
with a nontrivial operation $\rhd$, and $\Phi$ is not
a Lie group homomorphism. Its sub-adjacent operation is
\[
(x,y,z)\circ(u,v,w)=(x+u,y+v,w).
\]
\end{Example}

\begin{Example}\label{ex:engel-twisted-post-Lie}
Let $G=\mathbb{R}^4$ be the Engel Lie group, written in
exponential coordinates. For
$a=(x_1,x_2,x_3,x_4)$ and $b=(y_1,y_2,y_3,y_4)$,
define
$a\cdot b=(c_1,c_2,c_3,c_4)$,
where
\[
\begin{aligned}
c_1&=x_1+y_1,\qquad c_2=x_2+y_2,\qquad
c_3 =x_3+y_3+\frac12(x_1y_2-x_2y_1),\\
c_4&=x_4+y_4+\frac12(x_1y_3-x_3y_1) +\frac1{12}(x_1-y_1)(x_1y_2-x_2y_1).
\end{aligned}
\]
Its identity element is $1=(0,0,0,0)$, and
$(x_1,x_2,x_3,x_4)^{-1}=(-x_1,-x_2,-x_3,-x_4)$.
Its Lie algebra has a basis $E_1$, $E_2$, $E_3$, $E_4$ with
nonzero basic brackets
$[E_1,E_2]=E_3$,
$[E_1,E_3]=E_4$.
In these coordinates,
\[
\exp\left(\sum_{i=1}^4 x_iE_i\right)=(x_1,x_2,x_3,x_4).
\]

Define $\Phi\colon G\to G$ and $\rhd\colon G\times G\to G$ by
$\Phi(x_1,x_2,x_3,x_4)=(0,0,0,x_4)$ and
\[
(x_1,x_2,x_3,x_4)\rhd(y_1,y_2,y_3,y_4)
=
\left(
e^{x_4}y_1,
e^{-2x_4}y_2,
e^{-x_4}y_3,
y_4
\right).
\]
Then $(G,\cdot,\rhd,\Phi)$ is a twisted post-Lie group
with a nontrivial operation $\rhd$, and $\Phi$ is not
a Lie group homomorphism. Its sub-adjacent operation is
\[
(x_1,x_2,x_3,x_4)\circ(y_1,y_2,y_3,y_4)
=
\left(
e^{x_4}y_1,\,
e^{-2x_4}y_2,\,
e^{-x_4}y_3,\,
x_4+y_4
\right).
\]
The exponential compatibility holds globally:
\[
\Phi(\exp X)=\exp(\varphi(X)),
\qquad X\in\mathfrak{g},
\]
where $\varphi=\Phi_{\ast 1}$ is given by
\[
\varphi\left(\sum_{i=1}^4 x_iE_i\right)=x_4E_4.
\]
\end{Example}

\begin{Definition}
	Let $(G,\cdot_{G},\rhd_{G},\Phi_{G})$ and $(H,\cdot_{H},\rhd_{H},\Phi_H)$ be twisted post-Lie groups.
	A map $F\colon (G,\cdot_{G},\rhd_{G},\Phi_{G})\to (H,\cdot_{H},\rhd_{H},\Phi_H)$ is called a twisted post-Lie group homomorphism if $F$ is both a smooth map and a twisted post-group homomorphism.
	
\end{Definition}
Next, we show that the differentiation of a twisted post-Lie group is a twisted post-Lie algebra.
\begin{Theorem}\label{TPLL}
	Let $(G,\cdot,\rhd,\Phi)$ be a twisted post-Lie group. Then $(\mathfrak{g},[\cdot,\cdot],\triangleright,\varphi)$ is a twisted post-Lie algebra, where $\varphi=\Phi_{*1}$. Moreover, let $F\colon (G,\cdot_G,\rhd_G,\Phi_G)\to (H,\cdot_H,\rhd_H,\Phi_H)$ be a twisted post-Lie group homomorphism. Then $f=F_{*1_{G}}\colon (\mathfrak{g},[\cdot,\cdot]_{\mathfrak{g}},\triangleright_\mathfrak{g},\varphi_\mathfrak{g})\to (\mathfrak{h},[\cdot,\cdot]_{\mathfrak{h}},\triangleright_\mathfrak{h},\varphi_\mathfrak{h}) $ is a twisted post-Lie algebra homomorphism.
\end{Theorem}
\begin{proof}Since $L_1^\rhd$ is a group automorphism, $1\rhd 1=1$. Hence
$1\circ1=\Phi(1)\cdot(1\rhd1)=1$. Therefore, \eqref{TPG2} gives
$L_1^\rhd=(L_1^\rhd)^2$.
Since $L_1^\rhd$ is invertible, $L_1^\rhd=\operatorname{id}_G$.

Define
$\rho\colon  G\longrightarrow\operatorname{Aut}(\mathfrak g)$,
$\rho(a)=(L_a^\rhd)_{*1}$,
and let
$\alpha=(d\rho)_1\colon
\mathfrak g\longrightarrow\operatorname{Der}(\mathfrak g)$.
Then
$\alpha(x)(y)=x\triangleright y$.
Since $\alpha(x)\in\operatorname{Der}(\mathfrak g)$, identity~\eqref{TPL1} follows.

By Lemma~\ref{IDM2}, $\rho\circ\Phi=\rho$. Differentiating this
identity at $1$ gives $\alpha\circ\varphi=\alpha$, which proves~\eqref{TPL4}. Moreover, differentiating $\Phi^2=\Phi$ gives
$\varphi^2=\varphi$, and hence \eqref{TPL5} holds.

We now prove \eqref{TPL3}.
Since $L_a^\rhd$ is a Lie group automorphism, we have
$a\rhd\exp(sy)=\exp\bigl(s\rho(a)y\bigr)$.
Moreover,
$\rho(\exp(tx))y=y+t(x\triangleright y)+\cdots$.

Choose an exponential neighbourhood $V$ of $1$ in $G$ and
write $\log\colon V\to\mathfrak g$ for the corresponding local
logarithm. By restricting $t$ and $s$ if necessary, all the
expressions below lie in $V$. Using the exponential compatibility of $\Phi$ and the
Baker--Campbell--Hausdorff formula, we obtain
\[
\begin{aligned}
\left.\frac{\d}{\d t} \right|_{t=0}\left.\frac{\d}{\d s} \right|_{s=0}\log\bigl(\exp(tx)\circ\exp(sy)\bigr)=
\frac12[\varphi(x),y]+x\triangleright y.
\end{aligned}
\]
It follows that
\begin{align}
\left.\frac{\d}{\d t} \right|_{t=0}\left.\frac{\d}{\d s} \right|_{s=0}\left(
\log\Phi\bigl(\exp(tx)\circ\exp(sy)\bigr)
\right)&=
\varphi\left(\frac12[\varphi(x),y]+x\triangleright y\right),\nonumber\\
\left.\frac{\d}{\d t} \right|_{t=0}\left.\frac{\d}{\d s} \right|_{s=0}\left(
\log\Phi\bigl(\exp(sy)\circ\exp(tx)\bigr)
\right)&=\varphi\left(\frac12[\varphi(y),x]+y\triangleright x\right).\label{eq LE}
\end{align}

On the other hand,
\[
\begin{aligned}
\left.\frac{\d}{\d t} \right|_{t=0}\left.\frac{\d}{\d s} \right|_{s=0}\log\left(
\exp(tx)\circ\Phi(\exp(sy))
\right)=
\frac12[\varphi(x),\varphi(y)]
+x\triangleright\varphi(y)
,
\end{aligned}
\]
and
\[
\begin{aligned}
\left.\frac{\d}{\d t} \right|_{t=0}\left.\frac{\d}{\d s} \right|_{s=0}\log\left(
\exp(sy)\circ\Phi(\exp(tx))
\right)=
\frac12[\varphi(y),\varphi(x)]
+y\triangleright\varphi(x)
.
\end{aligned}
\]

By \eqref{TPG3},
$\Phi(a\circ b)=a\circ\Phi(b)$.
Applying the local logarithm to both sides of \eqref{TPG3} and
comparing the coefficients of $ts$ gives
\[
\varphi\left(\frac12[\varphi(x),y]+x\triangleright y\right)
=
\frac12[\varphi(x),\varphi(y)]
+x\triangleright\varphi(y)
\]
and
\[
\varphi\left(\frac12[\varphi(y),x]+y\triangleright x\right)
=
\frac12[\varphi(y),\varphi(x)]
+y\triangleright\varphi(x).
\]
Subtracting these identities yields
\[
\begin{aligned}
\varphi\left(
\frac12[\varphi(x),y]
+\frac12[x,\varphi(y)]
+x\triangleright y-y\triangleright x
\right)=
[\varphi(x),\varphi(y)]
+x\triangleright\varphi(y)
-y\triangleright\varphi(x).
\end{aligned}
\]
This proves \eqref{TPL3}.
Finally, we show \eqref{TPL2} holds.
For any $x,y,z\in \mathfrak{g}$, by \eqref{TPG2} and \eqref{FRT}, we~have
\begin{align*}
 x\triangleright(y\triangleright z)-y\triangleright (x\triangleright z)
	={}&\left.\frac{\d}{\d t} \right|_{t=0}\left.\frac{\d}{\d s} \right|_{s=0}\left.\frac{\d}{\d r} \right|_{r=0} L^{\rhd}_{\exp(tx)}L^{\rhd}_{\exp(sy)}(\exp(rz))\\
&-\left.\frac{\d}{\d t} \right|_{t=0}\left.\frac{\d}{\d s} \right|_{s=0}\left.\frac{\d}{\d r} \right|_{r=0}L^{\rhd}_{\exp(sy)}L^{\rhd}_{\exp(tx)}(\exp(rz))\\	
	={} &\left.\frac{\d}{\d t} \right|_{t=0}\left.\frac{\d}{\d s} \right|_{s=0}\left.\frac{\d}{\d r} \right|_{r=0}L^{\rhd}_{\exp(tx)\circ\exp(sy)}(\exp(rz))\\
&-\left.\frac{\d}{\d t} \right|_{t=0}\left.\frac{\d}{\d s} \right|_{s=0}\left.\frac{\d}{\d r} \right|_{r=0}L^{\rhd}_{\exp(sy)\circ\exp(tx)}(\exp(rz))\\	
	={}&\left.\frac{\d}{\d t} \right|_{t=0}\left.\frac{\d}{\d s} \right|_{s=0}\left.\frac{\d}{\d r} \right|_{r=0}L^{\rhd}_{\Phi(\exp(tx)\circ\exp(sy))}(\exp(rz))\\
&-\left.\frac{\d}{\d t} \right|_{t=0}\left.\frac{\d}{\d s} \right|_{s=0}\left.\frac{\d}{\d r} \right|_{r=0}L^{\rhd}_{\Phi(\exp(sy)\circ\exp(tx))}(\exp(rz)).
\end{align*}

Then it follows from \eqref{eq LE} that
\[
	x\triangleright(y\triangleright z)
	-y\triangleright(x\triangleright z)
	=
	\left(
	[\varphi(x),\varphi(y)]
	+x\triangleright\varphi(y)
	-y\triangleright\varphi(x)
	\right)\triangleright z.
\]
This proves \eqref{TPL2}.
Since $F\colon (G,\cdot_G)\to(H,\cdot_H)$ is a Lie group homomorphism, its differential
$f=F_{*1_G}$ is a Lie algebra homomorphism. For any $x,y\in\mathfrak{g}$, we~have
\begin{align*}
	f(x\triangleright_{\mathfrak{g}} y)&=f\left(\left.\frac{\d}{\d t} \right|_{t=0}\left.\frac{\d}{\d s} \right|_{s=0} \exp(tx)\rhd_G \exp(sy)\right) \\
	 & =\left.\frac{\d}{\d t} \right|_{t=0}\left.\frac{\d}{\d s} \right|_{s=0}F( \exp(tx)\rhd_G \exp(sy))\\
	&=\left.\frac{\d}{\d t} \right|_{t=0}\left.\frac{\d}{\d s} \right|_{s=0}F( \exp(tx))\rhd_H F(\exp(sy)).
\end{align*}
Then we have
\begin{align*}
	f(x\triangleright_{\mathfrak{g}} y)= \left.\frac{\d}{\d t} \right|_{t=0}\left.\frac{\d}{\d s} \right|_{s=0} \exp(tf(x))\rhd_H \exp(sf(y))
	= f(x)\triangleright_{\mathfrak{h}} f(y).
\end{align*}
And it follows from \eqref{TPHM3} that $f\circ \varphi_\mathfrak{g}=\varphi_{\mathfrak{h}}\circ f$.  This proves that $f$ is a homomorphism of twisted post-Lie algebras.\end{proof}

\begin{Lemma}\label{DAG}
	Let $(G,\cdot,\rhd,\Phi)$ be a twisted post-Lie group with the
	sub-adjacent operation $\circ$. Then $\Phi(G)$ is a Lie group
	with respect to $\circ$. Moreover, we have
	\[
	\left.\frac{\d}{\d t}\right|_{t=0}
	\Phi(\exp(tx))^{\dag}
	=
	-\varphi(x)
	=
	\left.\frac{\d}{\d t}\right|_{t=0}
	\Phi(\exp(tx))^{-1},
	\qquad \forall x\in\mathfrak g.
	\]
\end{Lemma}

\begin{proof}
By Lemma~\ref{IDM1}, we have $\Phi(1)=1$. Since $\Phi^2=\Phi$,
the map $\Phi$ is a smooth retraction onto $\Phi(G)$. Hence
$\Phi(G)$ is an embedded submanifold of $G$, and
$T_1\Phi(G)=\operatorname{Im}(\Phi_{*1}) =\varphi(\mathfrak g)$.
Since the multiplication $\circ$ and the inverse map
$\dagger$ restrict to smooth maps on $\Phi(G)$, Proposition~\ref{SG} implies that $(\Phi(G),\circ)$ is a Lie group with
identity~$1$.
Since $\cdot$, $\Phi$ and $\rhd$ are smooth, we obtain that $\circ$
is also smooth. By Lemma~\ref{IDM1}, we have $\Phi(1)=1$. It follows that $e_1=1$. Then by Proposition
\ref{SG} and Lemma~\ref{Phi}, we conclude that $(\Phi(G),\circ)$ is a Lie
group with identity $1$.

Finally, by the Leibniz rule and Lemma~\ref{Phi}, we have
\begin{align*}
\begin{split}
	0
	&=
	\left.\frac{\d}{\d t}\right|_{t=0}
	\bigl(
	\Phi(\exp(tx))^{\dag}\circ\Phi(\exp(tx))
	\bigr)
	 =
	\left.\frac{\d}{\d t}\right|_{t=0}
	\Phi(\exp(tx))^{\dag}
	+
	\left.\frac{\d}{\d t}\right|_{t=0}
	\Phi(\exp(tx))\\
	&=
	\left.\frac{\d}{\d t}\right|_{t=0}
	\Phi(\exp(tx))^{\dag}
	+
	\left.\frac{\d}{\d t}\right|_{t=0}
	\exp(t\varphi(x))
	 =
	\left.\frac{\d}{\d t}\right|_{t=0}
	\Phi(\exp(tx))^{\dag}
	+\varphi(x).
\end{split}
\end{align*}
Thus
\[
\left.\frac{\d}{\d t}\right|_{t=0}
\Phi(\exp(tx))^{\dag}
=
-\varphi(x).
\]
Similarly, differentiating
$\Phi(\exp(tx))^{-1}\cdot\Phi(\exp(tx))=1$
at $t=0$ gives
\[
\left.\frac{\d}{\d t}\right|_{t=0}
\Phi(\exp(tx))^{-1}
=
-\varphi(x).\tag*{\qed}
\]\renewcommand{\qed}{}
\end{proof}

Finally we prove that the Lie algebra of $(\Phi(G),\circ)$ is just $(\varphi(\mathfrak{g}),[\cdot,\cdot]_{\triangleright})$. This generalises \cite[Proposition 4.7]{BM}.
\begin{Proposition}
	Let $(G,\cdot,\rhd,\Phi)$ be a twisted post-Lie group. Then the Lie algebra of the Lie group $(\Phi(G),\circ)$ is $(\varphi(\mathfrak{g}),[\cdot,\cdot]_{\triangleright})$.
\end{Proposition}

\begin{proof}
For any $a,b\in G$, using Proposition~\ref{SG}, Lemmas~\ref{Phi} and~\ref{IDM2}, we have
\begin{align*}
	\Phi(a)\circ \Phi(b)\circ \Phi(a)^{\dag}&=\Phi(a)\circ \bigl(\Phi(b)\circ \Phi(a)^{\dag}\bigr)
	 =\Phi^{2}(a)\cdot\bigl(\Phi(a)\rhd\bigl(\Phi(b)\circ \Phi(a)^{\dag}\bigr)\bigr)\\
	&=\Phi(a)\cdot\bigl(\Phi(a)\rhd\bigl(\Phi(b)\cdot\bigl(\Phi(b)\rhd\Phi(a)^{\dag}\bigr)\bigr)\bigr).
\end{align*}

Denote the Lie bracket of the Lie group $(\Phi(G),\circ)$ by $[[\cdot,\cdot]]$. By the definition, we have
\begin{align*}
	[[\varphi(x),\varphi(y)]]&=\left.\frac{\d}{\d t}\right|_{t=0}\left.\frac{\d}{\d s}\right|_{s=0} \exp(t\varphi(x))\circ\exp(s\varphi(y))\circ\exp(t\varphi(x))^{\dag}\\
	&=\left.\frac{\d}{\d t}\right|_{t=0}\left.\frac{\d}{\d s}\right|_{s=0}\Phi(\exp(tx))\circ\Phi(\exp(sy))\circ\Phi(\exp(tx))^{\dag}\\
	&=\left.\frac{\d}{\d t}\right|_{t=0}\left.\frac{\d}{\d s}\right|_{s=0}\Phi(\exp(tx))\\
	&\quad\times\bigl(\Phi(\exp(tx))\rhd\bigl(\Phi(\exp(sy))\cdot\bigl(\Phi(\exp(sy))\rhd\Phi(\exp(tx))^{\dag}\bigr)\bigr)\bigr).
\end{align*}
Applying the differential of
the group multiplication, we obtain
\begin{align*}
[[\varphi(x),\varphi(y)]]
={}&
\left.\frac{\d}{\d t}\right|_{t=0}
\left.\frac{\d}{\d s}\right|_{s=0}
\Phi(\exp(tx))\\
&\quad{}\times
\bigl(\Phi(\exp(tx))\rhd\Phi(\exp(sy))\bigr)\cdot
\bigl(
\Phi(\exp(tx))\rhd
\Phi(\exp(tx))^\dagger
\bigr)\\
&+
\left.\frac{\d}{\d t}\right|_{t=0}
\left.\frac{\d}{\d s}\right|_{s=0}
\Phi(\exp(tx))\\
&\quad{}\times
\bigl(
\Phi(\exp(tx))\rhd
\bigl(
\Phi(\exp(sy))\rhd
\Phi(\exp(tx))^\dagger
\bigr)
\bigr).
\end{align*}
Since
\[
\Phi(\exp(tx))\rhd\Phi(\exp(tx))^\dagger
=
\Phi(\exp(tx))^{-1},
\]
the first term becomes
\begin{align*}
&\left.\frac{\d}{\d t}\right|_{t=0}
\left.\frac{\d}{\d s}\right|_{s=0}
\Phi(\exp(tx))\cdot
\bigl(\Phi(\exp(tx))\rhd\Phi(\exp(sy))\bigr)
\cdot\Phi(\exp(tx))^{-1}\\
&\qquad =
\left.\frac{\d}{\d t}\right|_{t=0}
\operatorname{Ad}^{\cdot}_{\Phi(\exp(tx))}
\bigl(
\bigl(L^\rhd_{\Phi(\exp(tx))}\bigr)_{*1}
\varphi(y)
\bigr)=
[\varphi(x),\varphi(y)]
+\varphi(x)\triangleright\varphi(y).
\end{align*}

For the second term, applying the chain rule with respect to $s$ and using Lemma~\ref{DAG}, we~obtain
\begin{align*}
&\left.\frac{\d}{\d t}\right|_{t=0}
\left.\frac{\d}{\d s}\right|_{s=0}
\Phi(\exp(tx))\cdot
\bigl(
\Phi(\exp(tx))\rhd
\bigl(
\Phi(\exp(sy))\rhd
\Phi(\exp(tx))^\dagger
\bigr)
\bigr)\\
&\qquad =
\left.\frac{\d}{\d s}\right|_{s=0}
\left.\frac{\d}{\d t}\right|_{t=0}
\Phi(\exp(sy))\rhd
\Phi(\exp(tx))^\dagger
 =
-\varphi(y)\triangleright\varphi(x).
\end{align*}

Combining the above identities gives
\begin{align*}
[[\varphi(x),\varphi(y)]]
&=
[\varphi(x),\varphi(y)]
+\varphi(x)\triangleright\varphi(y)
-\varphi(y)\triangleright\varphi(x)\\
&=
\bigl[\varphi^2(x),\varphi^2(y)\bigr]
+\varphi(x)\triangleright\varphi(y)
-\varphi(y)\triangleright\varphi(x)
=
[\varphi(x),\varphi(y)]_{\triangleright}.
\end{align*}
Thus $[[\cdot,\cdot]]=[\cdot,\cdot]_{\triangleright}$
on $\varphi(\mathfrak g)$.
\end{proof}

\section{(Weak) twisted post-Hopf algebras}\label{6}

In this section, first, we introduce the notion of (weak) twisted post-Hopf algebras. We show that the category of weak twisted post-Hopf algebras is isomorphic to the category of Hopf trusses proposed in \cite{TB2}. Then we show that every twisted post-Hopf algebra gives rise to another Hopf algebra, called the sub-adjacent Hopf algebra. Finally, we investigate the relationship between (weak) twisted post-groups and (weak) twisted post-Hopf algebras.

The definitions and results presented in this section were already
contained in the first version of the present preprint, which appeared
in 2023. Subsequently, Fern\'andez Vilaboa, Gonz\'alez Rodr\'{\i}guez and
Ramos P\'erez \cite{JF} extended the theory of weak twisted post-Hopf
algebras to arbitrary braided monoidal categories. They established
categorical relationships among weak twisted post-Hopf algebras, Hopf
trusses and weak twisted relative Rota--Baxter operators, and replaced
the global cocommutativity assumption by weaker compatibility
conditions related to the cocommutativity class of a Hopf algebra.
Below, we indicate more precisely how the results of this section are
related to these subsequent generalizations.

First, let us recall basic concepts on Hopf trusses.
In this section, we take a field of characteristic zero $\mathbb{F}$ as the base field.
\begin{Definition}
	Let $(H,\cdot,1,\Delta,\epsilon,S)$ be a cocommutative Hopf algebra. Let $\circ$ be a binary operation on $H$ such that $(H,\circ,\Delta,\epsilon)$ is a nonunital
	bialgebra. We say that $(H,\cdot,\circ,\Phi)$ is a Hopf truss if $\Phi\colon H\to H$ is a coalgebra homomorphism such that
	\begin{equation}\label{HTS}
		x\circ(y\cdot z)=(x_1\circ y)\cdot S(\Phi(x_2))\cdot(x_3\circ z), \qquad\forall x,y,z\in H.
	\end{equation}
	The operator $\Phi$ is called the cocycle of the Hopf truss $(H,\cdot,\circ,\Phi)$.
\end{Definition}
\begin{Definition}
	A Hopf truss morphism $F$ from $(H,\cdot_H,\circ_H,\Phi_H)$ to $(H',\cdot_{H'},\circ_{H'},\Phi_{H'})$ is a Hopf algebra morphism such that
	\begin{gather}\label{HT1}
		F(x\circ_H y)=F(x)\circ_{H'} F(y),
\\ \label{HT2}
		F(\Phi_H(x))=\Phi_{H'}(F(x)),
	\end{gather} for any $x,y\in H$.
\end{Definition}
\begin{Definition}
	Let $H$ be a cocommutative Hopf algebra with bilinear operator $\rhd\colon H\otimes H\to H$ and linear operator $\Phi\colon H\to H$. Define $\circ\colon H\otimes H\to H$ by
	\[x\circ y= \Phi(x_1)\cdot(x_2\rhd y),\quad\forall x,y\in H.\] $(H,\rhd,\Phi)$ is called a weak twisted post-Hopf algebra if $\rhd$ and $\Phi$ are coalgebra morphisms satisfying the following equalities:
\begin{subequations}
		\begin{gather}
			x\rhd (y\cdot z)=(x_1\rhd y)\cdot(x_2\rhd z),\label{TPHA1}\\
			(x\circ y)\rhd z=x\rhd (y\rhd z),\label{TPHA2}\\
			\Phi(x\circ y)=x\circ \Phi(y),\label{TPHA3}
		\end{gather}
	\end{subequations}
for any $x,y,z\in H$.
\end{Definition}
Moreover, if $\Phi(1)=1$ and the left multiplication $L^{\rhd}\colon H\to \End(H)$ given by
$L^{\rhd}_x(y)=x\rhd y$, $\forall y\in H$,  is invertible in $\Hom(H,\End(H))$ with respect to the convolution, that is, there is a~unique linear operator $T^{\rhd}\colon  H \to \End(H)$, such that for any $x\in H$,
\[L^{\rhd}_{x_1}\circ T_{x_2}^{\rhd}=T_{x_1}^{\rhd}\circ L^{\rhd}_{x_2}=\epsilon(x)\id_H,\] then $(H,\rhd,\Phi)$ is called a twisted post-Hopf algebra. The operator $\circ$ is called the sub-adjacent operation of the (weak) twisted post-Hopf algebra.

\begin{Definition}
	A (weak) twisted post-Hopf algebra homomorphism from (weak) twisted post-Hopf algebra $(H,\rhd,\Phi)$ to (weak) twisted post-Hopf algebra $(H',\rhd',\Phi')$ is a Hopf algebra morphism from $H$ to $H'$ such that
	\begin{gather}\label{TPH1}
		F(x\rhd y)=F(x)\rhd' F(y),\qquad\forall x,y\in H,\\
\label{TPH2}
		F(\Phi(x))=\Phi'(F(x)),\qquad\forall x\in H.
	\end{gather}
\end{Definition}

First, we show that every weak twisted post-Hopf algebra gives rise to a Hopf truss.
\begin{Proposition}\label{TPHT}
	Let $(H,\rhd,\Phi)$ be a weak twisted post-Hopf algebra. Define the operator $\circ\colon H\otimes H\to H$ by
$x\circ y=\Phi(x_1)\cdot(x_2\rhd y)$, $\forall x,y\in H$. Then $(H,\cdot,\circ,\Phi)$ is a Hopf truss. Moreover, let $F$ be a weak twisted post-Hopf algebra homomorphism from $(H,\rhd_H,\Phi_H)$ to $(H',\rhd_{H'},\Phi_{H'})$. Then~$F$ induces a Hopf truss morphism from $(H,\cdot_H,\circ_H,\Phi_H)$ to $(H',\cdot_{H'},\circ_{H'},\Phi_{H'})$.
\end{Proposition}
\begin{proof}
First, we show that $\circ$ is compatible with $\Delta$.
For any $x,y\in H$, we have
\begin{align*}
	\Delta(x\circ y)&=\Delta(\Phi(x_1)\cdot(x_2\rhd y))=(\Phi(x_1)\cdot(x_3\rhd y_1))\otimes(\Phi(x_2)\cdot(x_4\rhd y_2))\\
	&=(\Phi(x_1)\cdot(x_2\rhd y_1))\otimes\Phi(x_3)\cdot(x_4\rhd y_2)=(x_1\circ y_1)\otimes (x_2\circ y_2).
\end{align*}
This proves that $\circ$ is compatible with $\Delta$. Then we show that $\circ$ is compatible with $\epsilon$. We have
\[\epsilon(x\circ y)=\epsilon(\Phi(x_1)(x_2\rhd y))=\epsilon(\Phi(x_1))\epsilon(x_2)\epsilon(y)=\epsilon(x_1)\epsilon(x_2)\epsilon(y)=\epsilon(x)\epsilon(y).\] Next we show $(H,\circ)$ is a nonunital algebra.
For any $x,y,z\in H$, by the definition, we have
\begin{align*}
	x\circ (y\circ z)&=\Phi(x_1)\cdot(x_2\rhd(y\circ z))=\Phi(x_1)\cdot(x_2\rhd(\Phi(y_1)\cdot (y_2\rhd z)))\\
\overset{\eqref{TPHA1}}&{=}\Phi(x_1)\cdot(x_2\rhd \Phi(y_1))\cdot(x_3\rhd(y_2\rhd z))=(x_1\circ (\Phi(y_1)))\cdot(x_2\rhd(y_2\rhd z))\\
\overset{\eqref{TPHA3}}&{=}\Phi(x_1\circ y_1)\cdot((x_2\circ y_2)\rhd z)=(x\circ y)\circ z.
\end{align*}

Then we show \eqref{HTS} holds. By the definition, we have
\begin{align*}
	(x_1\circ y)\cdot S(\Phi(x_2))\cdot(x_3\circ z)&=\Phi(x_1)\cdot(x_2\rhd y)\cdot S(\Phi(x_3))\cdot\Phi(x_4)\cdot(x_5\rhd z)\\
	&=\Phi(x_1)\cdot(x_2\rhd y)\cdot(x_3\rhd z)\\
\overset{\eqref{TPHA1}}&{=}\Phi(x_1)\cdot(x_2\rhd (y\cdot z))=x\circ(y\cdot z).
\end{align*} This proves that \eqref{HTS} holds. Finally, if $F$ is a (weak) twisted post-Hopf algebra homomorphism from $(H,\cdot_H,\rhd_H,\Phi_H)$ to $(H',\cdot_{H'},\rhd_{H'},\Phi_{H'})$, then we have
\begin{align*}
	F(a\circ_H b)&=F(\Phi_H(a_1)\cdot_H(a_2\rhd_H b))=F(\Phi_H(a_1))\cdot_{H'}(F(a_2)\rhd_{H'} F(b))\\&=\Phi_{H'}(F(a_1))\cdot_{H'}(F(a_2)\rhd_{H'} F(b))
	=F(a)\circ_{H'}F(b).
\end{align*}
This implies \eqref{HT1} holds. And it follows from \eqref{TPH2} that \eqref{HT2} holds. Therefore, $F$ induces a Hopf truss morphism from $(H,\cdot_H,\circ_H,\Phi_H)$ to $(H',\cdot_{H'},\circ_{H'},\Phi_{H'})$.
\end{proof}

Conversely, we show that every Hopf truss gives rise to a weak twisted post-Hopf algebra.
\begin{Proposition}\label{HTTP}
	Let $(H,\cdot,\circ,\Phi)$ be a Hopf truss. Define
	$\rhd\colon H\otimes H\to H$ by
$x\rhd y=S(\Phi(x_1))\cdot(x_2\circ y)$, $\forall x,y\in H$.
	Then $(H,\rhd,\Phi)$ is a weak twisted post-Hopf algebra. Moreover, if $F$ is a~Hopf truss morphism from $(H,\cdot_H,\circ_H,\Phi_H)$ to $(H',\cdot_{H'},\circ_{H'},\Phi_{H'})$, then $F$ induces a weak twisted post-Hopf algebra homomorphism from $(H,\rhd_H,\Phi_H)$ to $(H',\rhd_{H'},\Phi_{H'})$
\end{Proposition}
\begin{proof}
First, we show that $\rhd$ is a coalgebra morphism.
For any $x,y\in H$, we have
\begin{align*}
	\Delta(x\rhd y)&=\Delta(S(\Phi(x_1))\cdot(x_2\circ y))=\left(S(\Phi(x_2))\cdot(x_3\circ y_1)\right)\otimes (S(\Phi(x_1))\cdot(x_4\circ y_2))\\
	&=(S(\Phi(x_1))\cdot(x_2\circ y_1))\otimes( S(\Phi(x_3))\cdot(x_4\circ y_2))=(x_1\rhd y_1)\otimes (x_2\rhd y_2).
\end{align*}	
This proves that $\rhd$ is a coalgebra morphism. It follows from \cite[Theorem 6.5]{TB2} that $(H,\rhd,\Phi)$ is a weak twisted post-Hopf algebra.
If $F$ is a Hopf truss morphism from $(H,\cdot_H,\circ_H,\Phi_H)$ to $(H',\cdot_{H'},\circ_{H'},\Phi_{H'})$, it is easy to see that \eqref{TPH2} holds, it remains to show \eqref{TPH1}. For any $x,y\in H$, we have
\begin{align*}
	F(x\rhd_{H} y)&=F(S_H(\Phi_H(x_1))\cdot_H(x_2\circ_H y))=F(S_H(\Phi_H(x_1)))\cdot_{H'}(F(x_2)\circ_{H'}F(y))\\
	&=S_{H'}(\Phi_{H'}(F(x_1)))\cdot_{H'}(F(x_2)\circ_{H'}F(y))
	=F(x)\rhd_{H'} F(y).\tag*{\qed}
\end{align*}\renewcommand{\qed}{}
\end{proof}

Let $\mathcal{HT}$ denote the category of Hopf trusses and $\mathcal{WTH}$ the category of weak twisted post-Hopf algebras. We have the following theorem.

Here is the Hopf algebra version of Theorem~\ref{th4}.
\begin{Theorem}\label{TPHH}
	The two categories $\mathcal{HT}$ and $\mathcal{WTH}$ are isomorphic.
\end{Theorem}
\begin{proof} By Proposition~\ref{TPHT}, we can define a functor $\mathcal{F}\colon \mathcal{WTH}\to \mathcal{HT}$. By Proposition~\ref{HTTP}, we can define a functor $\mathcal{G}\colon \mathcal{HT}\to \mathcal{WTH}$. It is easy to see that $\mathcal{F}\circ \mathcal{G}$ and $\mathcal{G}\circ \mathcal{F}$ are both identity maps. Hence, the two categories are isomorphic.
\end{proof}

\begin{Remark}
The categorical isomorphism in Theorem~\ref{TPHH} was established in
the original 2023 version of the present preprint for cocommutative
Hopf algebras in the category of vector spaces over a field.
Subsequently, \cite[Theorem~3.11]{JF} extended this correspondence to
an arbitrary braided monoidal category. More precisely, the authors
obtained an isomorphism between suitable full subcategories of weak
twisted post-Hopf algebras and Hopf trusses whose objects satisfy the
compatibility conditions in \cite[(3.7) and (3.14)]{JF}. When the
underlying monoidal category is symmetric, these conditions are
expressed in terms of the cocommutativity class and are weaker than
requiring the underlying Hopf algebra to be globally cocommutative.
\end{Remark}

Assume $(H,\rhd,\Phi)$ is a weak twisted post-Hopf algebra.
By the proof of Proposition~\ref{TPHT}, $(H,\circ,\Delta,\epsilon)$~is a nonunital bialgebra. Moreover, we have the following lemma.

\begin{Lemma}\label{YL}
	Let $(H,\rhd,\Phi)$ be a twisted post-Hopf algebra. Then for any $x,y\in H$, we have
	\begin{equation*}
		x\rhd 1=\epsilon(x)1,
\qquad
		1\rhd x=x,
\qquad
		S(x\rhd y)=x\rhd S(y).
	\end{equation*}
\end{Lemma}
\begin{proof}
Since $\rhd$ is a coalgebra morphism, for all $x,y\in H$ we have
\[
\Delta(x\rhd y)
=(x_1\rhd y_1)\otimes(x_2\rhd y_2),
\qquad
\epsilon(x\rhd y)=\epsilon(x)\epsilon(y).
\]

First, using the counit property, the antipode identity of $H$, and
\eqref{TPHA1}, we obtain
\begin{align*}
x\rhd 1
&=(x_1\rhd 1)\epsilon(x_2\rhd 1) =(x_1\rhd 1)\cdot(x_2\rhd 1)\cdot
 S(x_3\rhd 1)\\
&=(x_1\rhd 1)\cdot S(x_2\rhd 1)
=\epsilon(x\rhd 1)1
=\epsilon(x)1.
\end{align*}
Thus,
$x\rhd 1=\epsilon(x)1$.

Next, since $1$ is group-like, the convolution inverse identities
give
\[
L^\rhd_1\circ T^\rhd_1
=
T^\rhd_1\circ L^\rhd_1
=
\operatorname{id}_H.
\]
Hence $L^\rhd_1$ is invertible. Moreover, since $\Phi(1)=1$ and
$1\rhd1=1$, we have
\[
1\circ1=\Phi(1)\cdot(1\rhd1)=1.
\]
Therefore, by \eqref{TPHA2},
\[
(L^\rhd_1)^2(x)
=1\rhd(1\rhd x)
=(1\circ1)\rhd x
=1\rhd x
=L^\rhd_1(x).
\]
Thus $L^\rhd_1$ is both invertible and idempotent, and consequently
$L^\rhd_1=\operatorname{id}_H$. Hence
$1\rhd x=x$.

Finally, using the first identity and \eqref{TPHA1}, we obtain
\begin{align*}
S(x\rhd y)
&=S(x_1\rhd y_1)\epsilon(x_2)\epsilon(y_2) =S(x_1\rhd y_1)\cdot
(x_2\rhd(\epsilon(y_2)1))\\
&=S(x_1\rhd y_1)\cdot
(x_2\rhd(y_2\cdot S(y_3))) =S(x_1\rhd y_1)\cdot
 (x_2\rhd y_2)\cdot
 (x_3\rhd S(y_3))\\
&=\epsilon(x_1\rhd y_1)(x_2\rhd S(y_2))
=\epsilon(x_1)\epsilon(y_1)
 (x_2\rhd S(y_2))
=x\rhd S(y).
\end{align*}
Therefore,
$S(x\rhd y)=x\rhd S(y)$.
\end{proof}

\begin{Lemma}\label{IDMH} Let $(H,\rhd,\Phi)$ be a twisted post-Hopf algebra with the sub-adjacent operation $\circ$. Then the cocycle $\Phi$ is idempotent. Furthermore, for any $x,y\in H$, $\Phi(x)\rhd y=x\rhd y$.
\end{Lemma}
\begin{proof}
For any $x\in H$, we have
\[\Phi^{2}(x)=(x\circ 1)\circ 1=x\circ (1\circ1)=x\circ\Phi(1)=x\circ 1=\Phi(x).\] Then by \eqref{TPHA2} and Lemma~\ref{YL}, we have \[\Phi(x)\rhd y=(x\circ 1)\rhd y=x\rhd(1\rhd y)=x\rhd y.\tag*{\qed}
\]\renewcommand{\qed}{}
\end{proof}

\begin{Lemma}\label{CoalgebraT}
Let $(H,\rhd,\Phi)$ be a twisted post-Hopf algebra, and let
$T^\rhd\colon H\to\operatorname{End}(H)$ be the convolution inverse of
$L^\rhd$. Then, for all $x,y\in H$,
\begin{equation}\label{CoalgebraT1}
\Delta\bigl(T^\rhd_x(y)\bigr)
=
T^\rhd_{x_1}(y_1)\otimes T^\rhd_{x_2}(y_2),
\end{equation}
and
\begin{equation}\label{CoalgebraT2}
\epsilon\bigl(T^\rhd_x(y)\bigr)
=
\epsilon(x)\epsilon(y).
\end{equation}
\end{Lemma}

\begin{proof}
For every $x\in H$, define
$\mathcal{L}_x=L^\rhd_{x_1}\otimes L^\rhd_{x_2}$,
$\mathcal{T}_x=T^\rhd_{x_1}\otimes T^\rhd_{x_2}$.
Since $H$ is cocommutative, the convolution identities for
$L^\rhd$ and $T^\rhd$ imply
\[
\mathcal{L}_{x_1}\circ\mathcal{T}_{x_2}
=
\mathcal{T}_{x_1}\circ\mathcal{L}_{x_2}
=
\epsilon(x)\operatorname{id}_{H\otimes H}.
\]
Thus $\mathcal{T}$ is the convolution inverse of $\mathcal{L}$.

Since $\rhd$ is a coalgebra morphism, we have
$\Delta\circ L^\rhd_x=\mathcal{L}_x\circ\Delta$.
Applying $\Delta$ to
$L^\rhd_{x_1}\circ T^\rhd_{x_2}=\epsilon(x)\operatorname{id}_H$
gives
$
\mathcal{L}_{x_1}\circ\Delta\circ T^\rhd_{x_2}=\epsilon(x)\Delta$.
Convolving on the left with $\mathcal{T}$ and using the uniqueness of
the convolution inverse, we obtain
\[
\Delta\circ T^\rhd_x
=
\mathcal{T}_x\circ\Delta
=
\bigl(T^\rhd_{x_1}\otimes T^\rhd_{x_2}\bigr)\circ\Delta,
\]
which proves \eqref{CoalgebraT1}.

Moreover, since $\rhd$ is a coalgebra morphism,
$\epsilon (L^\rhd_x(y))=\epsilon(x)\epsilon(y)$.
Therefore,
\begin{align*}
\epsilon\bigl(T^\rhd_x(y)\bigr)
 =
\epsilon(x_1)\epsilon\bigl(T^\rhd_{x_2}(y)\bigr)
 =
\epsilon\bigl(
L^\rhd_{x_1}(T^\rhd_{x_2}(y))
\bigr)
 =
\epsilon(x)\epsilon(y).
\end{align*}
This proves \eqref{CoalgebraT2}.
\end{proof}

Now we give a Hopf algebra version of the sub-adjacent group.
\begin{Theorem}\label{SubAdjacentHopf}
	Let $(H,\rhd,\Phi)$ be a twisted post-Hopf algebra with the sub-adjacent operation~$\circ$. Then $(\Phi(H),\circ,\Delta,1,\epsilon)$ is a unital cocommutative bialgebra. Moreover, define $S_{\rhd}\colon \Phi(H)\to \Phi(H)$ by
\begin{equation}\label{rv1}
S_{\rhd}(\Phi(x))=\Phi(T^{\rhd}_{\Phi(x_1)}(S(\Phi(x_2)))),\qquad\forall x\in H.
 \end{equation}
Then $(\Phi(H),\circ,\Delta,1,\epsilon,S_{\rhd})$ is a cocommutative Hopf algebra.
\end{Theorem}
\begin{proof}
It follows from Proposition~\ref{TPHT} that $(H,\circ,\Delta,\epsilon)$ is a nonunital bialgebra.
By \eqref{TPHA3}, $\Phi(H)$~is compatible with $\circ$. And since $\Phi$ is a coalgebra morphism, $\Phi$ is compatible with $\Delta$, that is, for any $a\in H$, $\Delta(\Phi(a))=\Phi(a_1)\otimes \Phi(a_2)$.
Hence $(\Phi(H),\circ,\Delta,\epsilon)$ is a nonunital bialgebra.
Moreover, by Lemmas~\ref{IDMH} and~\ref{YL}, we have
$\Phi(x)\circ 1=\Phi^{2}(x)=\Phi(x)$, and $1\circ \Phi(x)=\Phi(1)\cdot(1\rhd \Phi(x))=\Phi(x)$ for any $x\in H$. This implies $(\Phi(H),\circ,\Delta,1,\epsilon)$ is a unital bialgebra.

Then we prove that $S_{\rhd}$ is well-defined. Define a linear map
\[
\Theta\colon \ H\otimes H\longrightarrow H,
\qquad
\Theta(u\otimes v)=T^\rhd_u(S(v)).
\]
Since $\Phi$ is a coalgebra morphism, for every $x\in H$ we have
\[
\Delta(\Phi(x))
=(\Phi\otimes\Phi)\Delta(x)
=\Phi(x_1)\otimes\Phi(x_2).
\]
Consequently,
\[
T^\rhd_{\Phi(x_1)}(S(\Phi(x_2)))
=\Theta\Delta(\Phi(x)).
\]
If $\Phi(x)=\Phi(y)$, then
$\Delta(\Phi(x))=\Delta(\Phi(y))$,
and hence
\[
T^\rhd_{\Phi(x_1)}(S(\Phi(x_2)))
=
T^\rhd_{\Phi(y_1)}(S(\Phi(y_2))).
\]
Thus \eqref{rv1} depends only on the element $\Phi(x)$ and not on the
choice of its representative~$x$.

Next, we verify that $S_{\rhd}$ is the antipode. It follows from Lemma~\ref{CoalgebraT} that
$x\longmapsto T^{\rhd}_{x_1}(S(x_2))$
is a coalgebra morphism. Since $\Phi$ is also a coalgebra morphism,
\eqref{rv1}~implies that $S_{\rhd}$ is a coalgebra morphism and
$\epsilon\circ S_{\rhd}=\epsilon$
on $\Phi(H)$.

By \eqref{TPHA3}, Lemma~\ref{IDMH}, and the convolution identity for
$T^{\rhd}$, we have
\begin{align}
\Phi(x_1)\circ S_{\rhd}(\Phi(x_2))
&=\Phi(x_1)\circ
 \Phi\bigl(T^{\rhd}_{\Phi(x_2)}(S(\Phi(x_3)))\bigr)
 =\Phi\bigl(
 \Phi(x_1)\circ
 T^{\rhd}_{\Phi(x_2)}(S(\Phi(x_3)))
 \bigr)\nonumber\\
&=\Phi\bigl(
 \Phi^2(x_1)\cdot
 \bigl(\Phi(x_2)\rhd
 T^{\rhd}_{\Phi(x_3)}(S(\Phi(x_4)))\bigr)
 \bigr)\label{C1}\\
&=\Phi\bigl(
 \Phi(x_1)\cdot
 \epsilon(\Phi(x_2))S(\Phi(x_3))
 \bigr)
 =\Phi\bigl(\Phi(x_1)\cdot S(\Phi(x_2))\bigr)
 =\epsilon(\Phi(x))1.\nonumber
\end{align}
Here, the second equality follows from \eqref{TPHA3}, namely,
$a\circ\Phi(b)=\Phi(a\circ b)$.

Note that
\begin{align}
S_{\rhd}^2(\Phi(x))
&=\epsilon(\Phi(x_1))S_{\rhd}^2(\Phi(x_2))
=(\Phi(x_1)\circ S_{\rhd}(\Phi(x_2)))
 \circ S_{\rhd}^2(\Phi(x_3))\nonumber\\
&=\Phi(x_1)\circ
 \bigl(S_{\rhd}(\Phi(x_2))
 \circ S_{\rhd}^2(\Phi(x_3))\bigr)
 =\Phi(x_1)\circ
 \epsilon(S_{\rhd}(\Phi(x_2)))1\nonumber\\
&=\Phi(x_1)\circ\epsilon(\Phi(x_2))1
 =\Phi(x).\label{C2}
\end{align}
Since $S_{\rhd}$ is a coalgebra morphism, applying \eqref{C1} to
$S_{\rhd}(\Phi(x))$ and using \eqref{C2}, we obtain
\begin{align*}
S_{\rhd}(\Phi(x_1))\circ\Phi(x_2)
 =
S_{\rhd}(\Phi(x_1))
\circ S_{\rhd}^{2}(\Phi(x_2))
 =
\epsilon(S_{\rhd}(\Phi(x)))1
 =
\epsilon(\Phi(x))1.
\end{align*}
Hence
$(\Phi(H),\circ,\Delta,1,\epsilon,S_{\rhd})$
is a cocommutative Hopf algebra.
\end{proof}

\begin{Remark}
The construction of the sub-adjacent Hopf algebra on
$\Phi(H)=\operatorname{Im}(\Phi)$
in Theorem~\ref{SubAdjacentHopf} already appeared in the original
2023 version of the present preprint. It was subsequently extended
in~\cite{JF} to the setting of braided monoidal categories.

More precisely, \cite[Lemma~3.13]{JF} proves that a unit-preserving
cocycle $\Phi_H$ is idempotent already for a weak twisted post-Hopf
algebra. If the underlying braided monoidal category admits split
idempotents, the idempotent $\Phi_H$ determines an image object
$I(\Phi_H)$ together with morphisms
\[
H\xrightarrow{\,p_H\,}I(\Phi_H)
 \xrightarrow{\,i_H\,}H,
\qquad
i_H\circ p_H=\Phi_H,
\qquad
p_H\circ i_H=\operatorname{id}_{I(\Phi_H)}.
\]
Under the compatibility condition \cite[(3.7)]{JF},
\cite[Theorem~3.21]{JF} equips $I(\Phi_H)$ with a bialgebra
structure. Under the additional assumptions of
\cite[Corollary~3.23]{JF}, this bialgebra becomes a cocommutative
Hopf algebra. This gives a braided
categorical extension of the image construction in Theorem~\ref{SubAdjacentHopf}.
\end{Remark}

The Hopf algebra $(\Phi(H),\circ,\Delta,1,\epsilon,S_{\rhd})$ given in the above theorem is called the sub-adjacent Hopf algebra of the twisted post-Hopf algebra $(H,\rhd,\Phi)$.

Finally, we show that the group-like elements of a~(weak) twisted post-Hopf algebra form a~(weak) twisted post-group.

\begin{Theorem}\label{TPHH1}
Let $(H,\rhd,\Phi)$ be a $($weak$)$ twisted post-Hopf algebra with the sub-adjacent operation $\circ$. Then its subset $G(H)$ of group-like elements has the structure of $($weak$)$ twisted post-group. Conversely, let $(G,\cdot,\rhd,\Phi)$ be a weak twisted post-group. Then the linear extensions of $\rhd$ and $\Phi$ make $\mathbb F[G]$ into a weak twisted post-Hopf algebra. If, moreover, $(G,\cdot,\rhd,\Phi)$ is a~twisted post-group and $\Phi(1_G)=1_G$, then $\mathbb F[G]$ is a~twisted post-Hopf algebra.
\end{Theorem}
\begin{proof}
Since $\Phi$ is a coalgebra morphism, for any $x\in G(H)$, we have
\[\Delta(\Phi(x))=\Phi(x_1)\otimes \Phi(x_2)=\Phi(x)\otimes\Phi(x).\] Hence we can define $\Phi\colon G(H)\to G(H)$.
Similarly, for any $x,y\in G(H)$, we have
\begin{align*}
	\Delta(x\rhd y)=(x_1\rhd y_1)\otimes (x_2\rhd y_2)=(x\rhd y)\otimes (x\rhd y).
\end{align*}
This implies we can define $\rhd\colon G(H)\times G(H)\to G(H)$ and define $\circ\colon G(H)\times G(H)\to G(H)$ given by $x\circ y=\Phi(x)\cdot(x\rhd y)$.
For any $x,y,z\in G(H)$, we have
\[x\rhd (y\cdot z)=(x_1\rhd y)\cdot(x_2\rhd z)=(x\rhd y)\cdot(x\rhd z).\] This implies
the left multiplication $L^{\rhd}_x\colon G(H)\to G(H)$ defined by
$L^{\rhd}_{x}(z)=x\rhd z$ for any $z\in G(H)$, is a group homomorphism. Moreover, if $(H,\rhd,\Phi)$ is a twisted post-Hopf algebra, then we have
\[L^{\rhd}_{x}T^{\rhd}_{x}=L^{\rhd}_{x_1}T^{\rhd}_{x_2}=\epsilon(x)\id=\id,\] and similarly we have $T^{\rhd}_{x}L^{\rhd}_{x}=\id$. This implies $L^{\rhd}_x$ is an automorphism.

For any $x,y,z\in G(H)$, by \eqref{TPHA2}, we have \eqref{TPG2}. And using \eqref{TPHA3}, we know that \eqref{TPG3} is satisfied. Hence $(G(H),\cdot,\rhd,\Phi)$ is a (weak) twisted post-group if $(H,\rhd,\Phi)$ is a (weak) twisted post-Hopf algebra.

Conversely, given a (weak) twisted post-group $(G,\rhd,\Phi)$, we can extend $\rhd$ and $\Phi$ to $\mathbb{F}[G]$ by
\[\sum_i \alpha_i a_i\rhd\sum_j\beta_j b_j=\sum_{i,j} \alpha_i \beta_j (a_i\rhd b_j),\qquad
\Phi(\sum_i \alpha_i a_i)=\sum_i \alpha_i \Phi(a_i),\]
where $\alpha_i,\beta_j\in \mathbb{F}$ and $a_i,b_j\in G$. One can readily verify $(\mathbb{F}[G],\rhd,\Phi)$ is a (weak) twisted post-Hopf algebra.
\end{proof}

\subsection*{Acknowledgments}
 The author would like to thank the anonymous referees for their helpful comments and suggestions, which have greatly improved the manuscript.
 This research is supported by Scientific Research Foundation for High-level Talents of Anhui University of Science and Technology (2024yjrc49).

\pdfbookmark[1]{References}{ref}
\LastPageEnding

\end{document}